\documentclass[a4paper,12pt]{amsart}
\usepackage{rotating,pdflscape,color,amssymb,hyperref,subfigure,psfrag,amsmath,eufrak,bbm,epsfig}
\usepackage{a4wide,amscd}
\usepackage{tikz} % doc : chercher PGF
 \usepackage[usenames,dvipsnames]{pstricks}
 \usepackage{pst-grad} % For gradients
 \usepackage{pst-plot} % For axes

\usepackage[small]{caption}

\newcommand{\ga}{\gamma}

\newcommand{\F}{\mc{F}}

\newcommand{\ld}{\ldots}
\newcommand{\beg}{\begin}
\newcommand{\en}{\end}

\newcommand{\trm}{\textrm}

\newcommand{\bgt}{\begin{itemize}}
\newcommand{\ent}{\end{itemize}}

\newcommand{\op}{\operatorname}
\newcommand{\eqre}{\eqref}
\newcommand{\re}{\ref}
\newcommand{\la}{\label}

\newcommand{\rfl}{\rfloor}
\newcommand{\lfl}{\lfloor}
\newcommand{\si}{\sigma}

\newcommand{\Var}{\operatorname{Var}}

\newcommand{\diag}{\operatorname{diag}}

\newcommand{\ds}{\displaystyle}

\newcommand{\p}{\mathbb{P}}

\newcommand{\Tr}{\operatorname{Tr}}

\newcommand{\ninf}{\underset{n\to\infty}{\longrightarrow}}

\newcommand{\E}{\op{\mathbb{E}}}

\newcommand{\R}{\mathbb{R}}
\newcommand{\C}{\mathbb{C}}

\newcommand{\ud}{\mathrm{d}}

\newcommand{\pro}{probability }

\newcommand{\f}{\frac}
\newcommand{\ff}{\frac{1}}
\newcommand{\lf}{\left}
\newcommand{\ri}{\right}

\newcommand{\st}{such that }

\newcommand{\lam}{\lambda}

\newcommand{\ti}{\times}

\newcommand{\var}{\Var}

\newcommand{\ste}{\, ;\, }
\newcommand{\mc}{\mathcal }

\newcommand{\eps}{\varepsilon}

\newcommand{\bck}{\backslash}
\newcommand{\al}{\alpha}

\newcommand{\ovl}{\overline}

\newcommand{\bbm}{\begin{bmatrix}}
\newcommand{\ebm}{\end{bmatrix}}
\newcommand{\bes}{\begin{equation*}}
\newcommand{\ees}{\end{equation*}}
\newcommand{\be}{\begin{equation}}
\newcommand{\ee}{\end{equation}}
\newcommand{\beqy}{\begin{eqnarray}}
\newcommand{\eeqy}{\end{eqnarray}}
\newcommand{\beq}{\begin{eqnarray*}}
\newcommand{\eeq}{\end{eqnarray*}}
\newcommand{\one}{\mathbbm{1}}
\newcommand{\lto}{\longrightarrow}

\newcommand{\ie}{\emph{i.e. }}

\newcommand{\bpm}{\begin{pmatrix}}
\newcommand{\epm}{\end{pmatrix}}

\newcommand{\cd}{\cdots}
\newcommand{\Lvy}{L\'evy }

\newcommand{\bpr}{\beg{proof}}
\newcommand{\epr}{\en{proof}}

\newcommand{\bet}{\beta}
\newcommand{\del}{\delta}
\newcommand{\Del}{\Delta}

\newcommand{\pa}{\partial}
\newcommand{\ra}{\rangle}

\newcommand{\ba}{\mathbf{a}}

\newcommand{\ka}{\kappa}

\newcommand{\tB}{\tilde{B}}

\newcommand{\bfe}{\mathbf{e}}

\newcommand{\cadlag}{c\`adl\`ag }
\newcommand{\sgn}{\op{sgn}}
\newcommand{\tz}{\tilde{z}}

\def\eq{\begin{eqnarray*}}
\def\qe{\end{eqnarray*}}
\def\eqa{\begin{eqnarray}}
\def\qea{\end{eqnarray}}

\newtheorem{Th}{Theorem}[section]

\newtheorem{propo}[Th]{Proposition}

\newtheorem{lem}[Th]{Lemma}

\theoremstyle{definition}

\newtheorem{Def}[Th]{Definition}

\long\def\symbolfootnote[#1]#2{\begingroup
\def\thefootnote{\fnsymbol{footnote}}\footnote[#1]{#2}\endgroup} 
\def\bC{{\mathbb C}}

\allowdisplaybreaks
\def\ra{{\rightarrow}}

\title{Central limit theorem  for eigenvectors of heavy tailed matrices}

\author[Florent Benaych-Georges, Alice Guionnet]{Florent Benaych-Georges, Alice Guionnet}

\thanks{FBG: florent.benaych-georges@parisdescartes.fr, MAP 5, UMR CNRS 8145 - Universit\'e Paris Descartes, 45 rue des Saints-P\`eres 75270 Paris Cedex~6, France.\\
AG: aguionne@ens-lyon.fr, CNRS \& \'Ecole Normale Sup\'erieure
de Lyon, Unit\'e de math\'ematiques pures et appliqu\'ees, 46 all\'ee
d'Italie, 69364 Lyon Cedex 07, France and MIT, Mathematics Department, 77 Massachusetts Av, Cambridge MA 02139-4307, USA.\\
Research supported by  Simons Foundation and the NSF award DMS-1307704.}


\begin{document}
\maketitle

\begin{abstract} We consider the eigenvectors  of symmetric matrices with independent heavy tailed entries, such  as matrices with 
entries in the domain of attraction of $\alpha$-stable laws,  or     adjacency matrices of  Erd\"os-R\'enyi graphs.
We denote by $U=[u_{ij}]$ the eigenvectors matrix (corresponding to increasing eigenvalues) and prove that the bivariate process
 $$B^n_{s,t}:=\ff{\sqrt{n}}\sum_{\substack{1\le i\le ns\\ 1\le j\le nt}}(|u_{ij}|^2 -\ff{n})\qquad\trm{($0\le s,t \le 1$)}, $$
   converges in law to  a non trivial  Gaussian process. An interesting part of this result is the $\ff{\sqrt{n}}$ rescaling, proving that from this point of view, the eigenvectors  matrix $U$ behaves more like a  permutation %(hence heavy tailed) 
 matrix (as it was proved in  \cite{chapuy} that for $U$ a permutation matrix, $\ff{\sqrt{n}}$ is the right scaling) than like a Haar-distributed orthogonal or unitary matrix (as it was proved in  \cite{cat-alain10,BeffaraDonatiRouault} that for $U$ such a  matrix,  the right scaling  is $1$).
\end{abstract}

%ArXiv abstract : 
%We consider the eigenvectors of symmetric matrices with independent heavy tailed entries, such as matrices with entries in the domain of attraction of $\alpha$-stable laws, or adjacencymatrices of Erdos-Renyi graphs. We denote by $U=[u_{ij}]$ the eigenvectors matrix (corresponding to increasing eigenvalues) and prove that the bivariate process $$B^n_{s,t}:=n^{-1/2}\sum_{1\le i\le ns, 1\le j\le nt}(|u_{ij}|^2 -n^{-1}),$$ indexed by $s,t\in [0,1]$, converges in law to  a non trivial Gaussian process. An interesting part of this result is the $n^{-1/2}$ rescaling, proving that from this point of view, the eigenvectors matrix $U$ behaves more like a  permutation matrix (as it was proved by Chapuy that for $U$ a permutation matrix, $n^{-1/2}$ is the right scaling) than like a Haar-distributed orthogonal or unitary matrix (as it was proved by Rouault and Donati-Martin that for $U$ such a matrix, the right scaling is $1$).

%\tableofcontents

\section{Introduction}

During the last decade, many breakthroughs were achieved in the study of random matrices belonging to the  GUE universality-class, that is Hermitian matrices with independent and equidistributed entries (modulo the symmetry constraint) with enough finite moments. 
The first key result about such matrices is due to Wigner \cite{wigner} in the fifties who showed 
  that the macroscopic behavior of
their eigenvalues is universal and asymptotically described by the semi-circle distribution. However, it took a long time to get more precise information on the local behavior of the eigenvalues, and for instance about the asymptotic distribution of their spacings. Even though local results were conjectured, for instance by Dyson and Mehta \cite{mehta}, it is only in the nineties that rigorous results were derived, such as the convergence of the  joint probability distribution  of eigenvalues in an interval of size of order $N^{-1}$ or the fluctuations of the largest eigenvalues, see \cite{TW1}.
Yet these results were restricted to Gaussian ensembles for which the joint law of the eigenvalues is known. Recently, these results were shown to be universal, that is to hold also for matrices with independent non Gaussian entries, provided they have enough finite moments \cite{U1,U2,U3,U4,TV}. Such a simple question as the convergence of the
law of a single spacing was   open, even in the GUE case,  until recently when it was solved by Tao \cite{Tao}. Once considering non Gaussian matrices,  it is natural to wonder about the behavior of the eigenvectors and whether they are delocalized (that is go to zero in $L^\infty$ norm)  as for GUE matrices. This was indeed shown by Erd\"os, Schlein and Yau \cite{Udeloc}.

Despite the numerous breakthroughs concerning random matrices belonging to the  GUE universality-class, not much is yet known about  other matrices. A famous example of
such a matrix is given by the adjacency matrix of an Erd\"os-R\'enyi graph. Its entries are independent (modulo the symmetry hypothesis) and equal to one with probability 
$p=p(N)$, zero otherwise.  If $pN$ goes to infinity fast enough, the matrix belongs to the GUE universality class \cite{Uerdos}. However if $pN$ converges to  a finite   non zero
constant, the matrix behaves quite differently, more like  a ``heavy tailed random matrix'', \ie a matrix  filled with independent entries
which have no finite second moment. Also in this case, it is known that, once properly normalized,  the empirical measure of the eigenvalues converges weakly  but the limit differs from the semi-circle distribution \cite{ZAK, BAGheavytails, BDG, FT, charles_alice, ACFTCL}.  Moreover, the fluctuations of the empirical measure could be studied \cite{ACFTCL, tirozzi, KKP96, KSV04}. It turns out that it fluctuates much more than in the case of matrices from the  GUE universality-class, as fluctuations are square root of the dimension bigger.
However, there is no result about the local fluctuations of the eigenvalues except in the case of matrices with entries in the domain of attraction of an $\alpha$-stable law in which case it was shown \cite{abap, sash}  that the largest eigenvalues are much bigger than the others, converge to a Poisson distribution and have localized eigenvectors.  About localization and delocalization of the eigenvectors, some models are  conjectured \cite{CB, slanina} to exhibit a phase transition; eigenvalues in a compact would have more delocalized eigenvectors than
outside this compact. Unfortunately, very little could be proved so far in this direction. Only the case where the entries are $\alpha$-stable random variables could be tackled \cite{charles_alice}; it was shown that for $\alpha>1$ the eigenvectors are delocalized whereas for $\alpha<1$ and large eigenvalues, a weak form of localization holds.

In this article, we study another type of properties of the eigenvectors of  a random matrix. Namely we consider the bivariate process
$$G^n_{s,t}:=\sum_{\substack{1\le i\le ns\\ 1\le j\le nt}}(|u_{ij}|^2 -\ff{n})\qquad\trm{($0\le s \le 1$, \quad $0\le t \le 1$)}\,,$$ where $U=[u_{ij}]$ is an   orthogonal    matrix whose columns are the eigenvectors of an Hermitian  random matrix $A=[a_{k\ell}]$. 
In the case where $A$ is a GUE matrix \cite{cat-alain10}, and then a more general matrix in the GUE universality-class \cite{FloEigenvectors}, it was shown that this process
converges in law towards a bivariate Brownian bridge (see also the closely related issues considered in \cite{BaoPanZhouUnivWigner}).   Here, we investigate the same process in the case where $A$ is a heavy tailed random matrix and show
that it fluctuates much more, namely it is $\ff{\sqrt{n}}G^n$ which converges in law. The limit  is a Gaussian process whose covariance depends 
on the model through the function $$\Phi(\lam)=\lim_{n\to\infty} n\E[e^{-i\lam a_{k\ell}}-1]\,.$$
 Of course if one considers heavy tailed variables in the domain of attraction of the same $\alpha$-stable law, and they suitably renormalize them, the function $\Phi$ and therefore the covariance will be the same. 
However, the covariance may vary when $\Phi$ does and it is not trivial in the sense that it does not vanish uniformly if $\Phi$ is not linear (a case which corresponds to light tails). 
In the related  setting 
when the $u_{ij}$'s are the entries of a uniformly distributed random permutation, \ie of a somehow sparse matrix, the process $\ff{\sqrt{n}}G^n$ also converges in law, towards the bivariate Brownian bridge   \cite{chapuy}. 

More precisely, we consider a real symmetric random $n\ti n$ matrix $A$ that can be either a Wigner matrix with exploding moments (which includes the adjacency matrix for Erd\"os-R\'enyi graphs) or a matrix with i.i.d. entries in the domain of attraction of a stable law (or more generally a matrix satisfying the hypotheses detailed in Hypothesis \re{Hyp:Model}). We then introduce an orthogonal    matrix $U=[u_{ij}]$ whose columns are the eigenvectors of $A$ so that we have 
$A=U\diag(\lam_1, \ld, \lam_n)U^*$. We then define the bivariate processes $$B^n_{s,t}:=\ff{\sqrt{n}}\sum_{\substack{1\le i\le ns\\ 1\le j\le nt}}(|u_{ij}|^2 -\ff{n})\qquad\trm{($0\le s \le 1$, \quad $0\le t \le 1$)}$$ and $$C^n_{s,\lam}:=\ff{\sqrt{n}}\sum_{\substack{1\le i\le ns\\ 1\le j\le n\ste \lam_j\le \lam}}(|u_{ij}|^2 -\ff{n})\qquad\trm{($0\le s \le 1$, \quad $\lam\in \R$)}$$ and prove, in Theorem \re{mainresult30513}, that both of these processes (with a little technical restriction on the domain of $B$) converge in law to (non trivial) Gaussian processes linked by the relation $$B_{s,F_{\mu_\Phi}(\lam)}=C_{s,\lam},$$ where $F_{\mu_\Phi}(\lam)=\mu_\Phi((-\infty,\lam])$ denotes the cumulative distribution function of the limit spectral law $\mu_\Phi$ of $A$,  \ie \be\la{BuSpFi}F_{\mu_\Phi}(\lam)=\lim_{n\to\infty} F_n(\lam), \quad \trm{ with } \quad F_n(\lam):=\ff{n}|\{i\ste \lam_i\le \lam\}|.\ee

The idea of the proof is the following one. We first notice that   for any $s\in [0,1]$,      the function      $\lam\mapsto C^n_{s,\lam}$ is the cumulative distribution function of the  random signed measure $ \nu_{s,n} $ on $\R$ defined by  \be\la{BSFDnu}\nu_{s, n}:= \ff{\sqrt{n}}\sum_{1\le i\le ns}\sum_{j=1}^n \lf(|u_{ij}|^2-\ff{n}\ri) \del_{\lam_j}\ee (\ie that  for any $\lam\in\R$, $\ds 
 C^n_{s,\lam}=\nu_{s,n}((-\infty, \lam])$).
 Then, we introduce the Cauchy transform $\ds X^n(s,z):=\int_{\lam\in \R}\f{\ud \nu_{s,n}(\lam)}{z-\lam}$ of $\nu_{s,n}$ and prove (Proposition \re{PropoCVST}) that  the process $(X^n(s,z))_{s,z}$ converges in law to a limit Gaussian process $(H_{s,z})$. This convergence is proved thanks to the classical CLT for martingales (Theorem \re{thconvmart} of the Appendix) together with the Schur complement formula  and  fixed points characterizations like the ones of the papers \cite{BAGheavytails,BDG,ACFTCL}. Then to deduce the convergence in law of the process $(C^{n}_{s,\lam})_{s, \lam}$, we use  the idea that the cumulative distribution function of a signed measure is entirely determined by its Cauchy transform. In fact, as   the measures $\nu_{s,n}$ of \eqre{BSFDnu} are random, things are slightly more complicated, and we need to prove a tightness lemma for the process    $(C^{n}_{s,\lam})_{s, \lam}$ (specifically Lemma \re{lemCtightness} of the Appendix, first applied to the process $(B^{n}_{s,t})$ and then transferred to 
 $(C^{n}_{s,\lam})$ by   Formula \eqre{89131} below). This lemma reduces the problem to the proof of the unicity of the possible limits for $(C^{n}_{s,\lam})$.  Then, we   use the formula  \bes\la{305132ABCDE} \int_{\lam\in \R}\f{ C^n_{s, \lam}}{(z-\lam)^2}\ud \lam  =-X^n(s,z)\ees and Lemma  \re{lemcauchytransfo} of the appendix to be able to claim that  $(C^{n}_{s,\lam})_{s, \lam}$ has a unique  limit  point. 
 The result proved for $(C^{n}_{s,\lam})$ can then be transferred to $(B^n_{s,t})$ thanks to  Formula \eqre{89131} below, where  $F_n(\lam)$, defined at \eqre{BuSpFi} above,   converges to the deterministic limit  $F_{\mu_\Phi}(\lam)$: 
  \be\la{89131}C^n_{s,\lam}=B^n_{s,F_n(\lam)}  .
\ee  

\beg{rmk}\la{2441416h}The objects introduced in the previous paragraph for the sketch of the proof enlighten the reason of the presence of the factor $\ff{\sqrt{n}}$ in the definitions of the processes $B^n_{s,t}$ and $C^n_{s,t}$ (recall that this factor does not appear in the corresponding formulas when Wigner instead of heavy-tailed matrices are concerned). 
Let $\mu_n$ denote the empirical spectral law of $A$ and let, for $i=1, \ld, n$,  $\mu_{n,\bfe_i}$ denote the empirical spectral law of $A$ according to the $i$th vector $\bfe_i$ for the canonical basis: that is, for any test function $f$, \be\la{2441411h31}\int f(x)\ud \mu_{n}(x)=\ff{n}\Tr(f(A))\qquad;\qquad \int f(x)\ud \mu_{n,\bfe_i}(x)=\sum_{j=1}^n|u_{i,j}|^2f(\lam_j)=(f(A))_{ii}.\ee
Then $$C^n_{s,\lam}=\ff{\sqrt{n}}\sum_{i=1}^{ns}(\mu_{n, \bfe_i}-\mu_n)((-\infty, \lam]),$$ so that $C^n_{s,\lam}$ is the centered version of a sum of  random variables $\mu_{n, \bfe_i}((-\infty, \lam])$ ($1\le i\le ns$). 
It has been proved in 
\cite{BCCAOP} that in the \Lvy case,  the random \pro measures $\mu_{n, \bfe_i}$ converge to i.i.d. copies of a nontrivial limiting random \pro measure (the spectral measure at the root of a suitable random weighted tree). This contrasts with case of
Wigner matrices, where concentration implies that the limiting measure of $\mu_{n, \bfe_i}$ is
deterministic, and  give a heuristic explanation, in the \Lvy case, of why one has to renormalize by $\sqrt{n}$ in $C^n_{s,\lam}$. Note however   that this explanation is not enough to prove that the variance of  $C^n_{s,\lam}$ does not explode nor vanish because $C^n_{s,\lam}$ is a sum of a large number of $\mu_{n, \bfe_i}((-\infty, \lam])$'s, that are correlated at finite $n$ (for example because the process vanishes on the boundary).
\en{rmk}

\noindent{\bf Organization of the paper.} The main results are stated in Section \re{sec:mainresults}. In  Section \re{secProof1}, we give a proof of Theorem  \re{mainresult30513}, based on Proposition \re{PropoCVST}, which is proved in Section \re{sec:PropoCVST}. Proposition \re{NonZeroCovariance}   is  proved   in Section \re{sec:NonZeroCovariance}. At last, some technical results are proved or recalled in the Appendix.

\noindent {\bf Notations.} For $u$, $v$ depending implicitly on $n$, we write $u\ll v$ when $u/v\lto 0$ as $n\to\infty$. For $x$ a random variable, $\var(x)$ denotes the \emph{variance} of $x$, \ie $\E[|x|^2]-|\E x|^2$. Power functions   are defined on $\C\bck\R_-$ via the standard   determination of the argument  on this set taking values in $(-\pi,\pi)$. The set $\C^+$ (resp. $\C^-$) denotes the open upper (resp. lower)  half plane and for any $z\in \C$, $\op{sgn}_{z}:=\op{sign}(\Im z)$. At last, bor any variable $x$, $\pa_x$ denotes $\f{\pa}{\pa x}$.

\section{Main results}\la{sec:mainresults}

Although technical, the model introduced in Hypothesis \ref{Hyp:Model} below has the advantage to be general enough to  contain several  models of interest. 
 \beg{hyp}\la{Hyp:Model}Let, for each $n\ge 1$, $A_n=[a_{ij}]$ be an $n\ti n$ real  symmetric   random matrix whose  sub-diagonal entries are some i.i.d. copies of a random variable $a$ (depending implicitly on $n$) such that:\\
$\bullet$ The random variable $a$ can be decomposed into $a=b+c$ \st  as $n\to\infty$, 
	  \beqy
  		\la{107130bis} &&\p(c\ne 0) \ll n^{-1}\\
		\la{107131bis} &&\var(b)\ll n^{-1/2}
	\eeqy
Moreover,  if the $b_i$'s are independent copies of $b$,
\begin{equation}\label{tyu}
\lim_{K\ra\infty}\lim_{n\ra\infty} \mathbb P\left(\sum_{i=1}^n (b_i-\mathbb E(b_i))^2\ge K\right)=0\,.\end{equation}
%A: ai ajoute cette condition car il n'est pas clair que les hop suffisentn
  $\bullet$ For any $\eps>0$ independent of $n$, the random variable $a$ can be decomposed into $a=b_\eps+c_\eps$ \st 
	\be\la{limceps}
		\limsup_{n\to\infty} n\, \p(c_\eps\ne 0)\le \eps
	\ee
for all $k\ge 1$, $n\E[(b_\eps-\E b_\eps)^{2k}]$ has a finite limit $C_{\eps,k}$ as $n\to\infty$. 

$\bullet$ For $\phi_n$ the function defined on  the closure $\ovl{\C^-}$ of $\C^- := \{ \lambda \in \C \ste \Im \lambda < 0\}$ by 
	\be\la{2071216h23} 
		\phi_n(\lambda) := \E \big[ \exp( -i \lambda a^2) \big], 
	\ee we have the convergence, uniform on  compact subsets of  $\ovl{\C^-}$, 
  \be\la{2071216h33}  n(\phi_n(\lam)-1)\;\lto\;\Phi(\lam),
\ee
 for a certain function $\Phi$ defined on $\ovl{\C^-}$.
 
 $\bullet$ The function $\Phi$ of \eqre{2071216h33}  admits the decomposition \be\la{hypcalcconv}
		\Phi(z)=\int_0^\infty g(y) e^{i\frac{y}{z}} \ud y
	\ee
	where $g(y)$ is a function  \st for some constants $K, \ga>-1, \ka\ge 0$, we have 
	\be\la{ConditionOng}
		|g(y)| \le K \one_{y\le 1} y^\ga +K \one_{y\ge 1} y^{\ka}, \qquad \forall y>0.
	\ee

$\bullet$ The function $\Phi$ of \eqre{2071216h33} also  either has the form \be\la{uniquenessassumption2}\Phi(x)=-\sigma (ix)^{\alpha/2} \ee or admits the (other) decomposition, for $x,y$ non zero:  \be\label{phias}
\Phi(x+y)=\iint_{(\R_+)^2}e^{i\frac{v}{x}+i\frac{v'}{y}} \ud\tau(v,v')+ \int_{\R^+} e^{i\frac{v}{x}} \ud\mu(v) + \int_{\R^+} e^{i\frac{v'}{y}} \ud\mu(v') \end{equation}	for some  complex measures $\tau,\mu$ on respectively $(\R^+)^2$  and $\R^+$ \st for all $b>0$, $\int e^{-bv}\ud|\mu|(v)$ is finite and for some constants $K>0$, $-1< \gamma \leq 0$ and $\ka \geq 0$, and		
	\eqa\label{eq:AssumpTau}
		 \frac{\ud|\tau|(v,v')}{\ud v \ud v'}\le K \big( v^\ga \one_{v\in ]0,1]} + v^{{\ka}} \one_{v\in ]1,\infty[}\big) \big( {v'}^{\ga} \one_{v'\in ]0,1]} + {v'}^{{\ka}} \one_{v'\in ]1,\infty[}\big).
	\qea
 \en{hyp}
 
 \beg{rmk}\la{22101318h17}When $\Phi$ satisfies \eqre{uniquenessassumption2} (e.g. for L\'evy matrices),
\eqref{phias} holds as well. Indeed,  for all $x,y\in\mathbb C^+$ (with a constant $C_\al$ that can change at every line),
\begin{eqnarray}\la{223131}
&&\Phi(x^{-1}+y^{-1})=C_\al (\frac{1}{x}+\frac{1}{y})^{\alpha/2}= C_\al \frac{1}{x^{\alpha/2}} \frac{1}{y^{\alpha/2}} (x+y)^{\alpha/2}\\
\nonumber &=&
 C_\al\int_0^\infty \ud w \int_0^\infty \ud w'\int_0^\infty \ud v
w^{\alpha/2-1} {w'}^{\alpha/2-1} v^{-\alpha/2-1} e^{ i w x+iw'y}(e^{iv(x+y)}-1)
\end{eqnarray}
(where we used the formula $\ds z^{\al/2} =C_\al\int_{t=0}^{+\infty}\f{e^{itz}-1}{t^{\al/2+1}}\ud t$ for any $z\in \C^+$ and $\al\in (0,2)$, which can be proved with the residues formula) so that \eqre{phias}   holds with $\mu= 0$ and $\tau(v,v')$ with density with respect to Lebesgue measure given by
	\beqy\label{eq:TauLevy}
	&&C_\al\int_0^{+\infty}  u^{-\alpha/2-1} \{(v-u)^{\alpha/2-1}(v'-u)^{\alpha/2-1}\one_{0\le u\le v\wedge v'}- v^{\alpha/2-1} {v'}^{\alpha/2-1}\}\ud u.
	\eeqy
Unfortunately, $\tau$ does not satisfy \eqref{eq:AssumpTau} as its density blows up at $v=v'$: we shall treat both case separately.\en{rmk}
 
 Examples of random matrices satisfying Hypothesis \re{Hyp:Model} are defined as follows.

\begin{Def}[Models of symmetric heavy  tailed matrices]~\label{defA}
Let $A=(a_{i,j})_{i,j=1,\ld,  n}$ be a random symmetric matrix with i.i.d. sub-diagonal entries.
\begin{enumerate}	\item[1.]
	We say that $A$ is a {\bf L\'evy matrix} of parameter $\alpha$ in $]0,2[$ when $A=X/a_n $ where   the entries $x_{ij}$ of $X$   have absolute values in the domain of attraction 
	of $\al$-stable distribution, more precisely 
	\be\la{ABP09exponent}\p\left(|x_{ij}|\ge u\right)=\frac{L(u)}{u^\alpha}\ee
	with a slowly varying function $L$, and
	$$a_n=\inf\{u: P\left(|x_{ij}|\ge u\right)\le\frac{1}{n}\}$$
($a_n=\tilde{L}(n)n^{1/\al}$, with $\tilde{L}(\cdot)$ a slowly varying function\footnote{A function $\tilde{L}$ is said to be  \emph{slowly varying} $\tilde{L}$ if for any fixed $\lam>0$, $\tilde{L}(\lam x)/\tilde{L}(x)\lto 1$ as $x\to\infty$.}).\\
	\item[2.]
		 We say that $A$ is a {\bf Wigner matrix with exploding moments} with parameter $(C_k)_{k\geq 1}$ whenever the entries of $A$ are centered, and  for any $k\geq 1$
	\be\la{1971214h}
		n\E\big[ (a_{ij})^{2k}\big]\ninf C_k>0.
	\ee    
 We assume that there exists a unique measure $m$  on $\mathbb R^+$ such that for all $k\ge 0$, \be\la{581320h30}C_{k+1}=\int x^k \ud m(x).\ee 
	\end{enumerate}
\end{Def}

The following proposition has been proved at Lemmas 1.3, 1.8 and 1.11 of \cite{ACFTCL}. 

\beg{propo}Both \Lvy matrices and Wigner matrices with exploding moments satisfy Hypothesis \re{Hyp:Model}:

$\bullet$ For \Lvy matrices, the function $\Phi$ of \eqre{2071216h33}  is given by formula  
 \be\la{exampleintroHTAD}\Phi(\lam)=-\sigma (i\lambda)^{\alpha/2}\ee
for some constant $\sigma> 0$,  the function $g$ of \eqre{hypcalcconv} is  $g(y)=C_\al y^{\f{\al}{2}-1}$,  with $C_\alpha=-\sigma i^{\alpha/2}$.

$\bullet$ For Wigner matrices with exploding moments, the function $\Phi$ of \eqre{2071216h33} is given by  \be\la{exampleWMWEMHTAD}\Phi(\lam)=\int \underbrace{\f{e^{-i\lam x}-1}{x}}_{:=-i\lam\trm{ for $x=0$}}\ud m(x),\ee for $m$ the measure of \eqre{581320h30},  the function $g$ of \eqre{hypcalcconv} is  
  \be\la{exampleWMWEMHTADLGI}g(y)=-\int_{\R^+} \underbrace{\f{J_1(2\sqrt{xy})}{\sqrt{xy}}}_{:=1\trm{ for }xy=0}\ud m(x),\ee for  $J_1$ the Bessel function of the first kind  defined by $\ds J_1 (s) = \frac s 2 \sum_{k \geq 0} \frac { (-s^2 / 4) ^k } {k ! (k+1)!}$,  and the measures  $\tau$ and $\mu$  of \eqre{phias} are  absolutely continuous with densities     \be\la{tauWMWEM}
 \f{\ud \tau(v,v')}{\ud v\ud v'}:=\int\f{J_1(2\sqrt{vx})J_1(2\sqrt{v'x})}{\sqrt{vv'}}\ud  m(x) \quad  \trm{ and }\quad \f{\ud \mu(v)}{\ud v}:=- \int\f{J_1(2\sqrt{vx})}{\sqrt{v}}\ud  m(x).\ee   %and   \be\la{muWMWEM}\f{\ud \mu(v)}{\ud v}:=- \int\f{J_1(2\sqrt{vx})}{\sqrt{v}}\ud  m(x).\ee 
\en{propo}

One can   easily see that our results also apply to complex Hermitian matrices: in this case, one only needs to require Hypothesis \re{Hyp:Model} to be satisfied by the absolute value of non diagonal entries and to have $a_{11}$ going to zero  as $N\to\infty$.

  A L\'evy matrix whose entries are truncated in an appropriate way is a Wigner matrix with exploding moments \cite{BAGheavytails, MAL122, ZAK}. The   recentered version\footnote{The recentering has in fact asymptotically  no effect on  the spectral measure $A$ as it is a rank  one   perturbation.} of the adjacency matrix of an Erd\"os-R\'enyi graph, \ie of a matrix $A$ \st   \be\la{exampleintroErdosReniymatrixAD}\trm{$A_{ij}=1$ with probability $p/n$ and $0$ with probability $1-p/n$},\ee is also an exploding moments Wigner matrix, with $\Phi(\lambda)=p(e^{-i\lambda}-1)$ (the measure $m$   is $p\del_1$). In this case the fluctuations were already studied in \cite{tirozzi}.

It has been proved in \cite{ACFTCL} (see also \cite{ZAK,tirozzi}) that under Hypothesis \re{Hyp:Model},  
 the    empirical spectral law \be\la{3051ter}\mu_n:=\ff{n}\sum_{j=1}^n \del_{\lam_j}\ee converges weakly in \pro to a deterministic \pro measure $\mu_\Phi$ that depends only on $\Phi$, \ie that for any continuous bounded function $f:\R\to \C$, we have the almost sure  convergence  \be\la{77132}\ff{n}\Tr f(A)=\ff{n}\sum_{j=1}^n f(\lam_j)\ninf \int f(x)\ud \mu_\Phi(x).\ee
  We introduce $F_{\mu_\Phi}(\lam):=\mu_\Phi((-\infty, \lam])$, cumulative distribution function of $\mu_\Phi$, and define the set $E_\Phi\subset[0,1]$   by \be\la{defE_Phi}E_\Phi:=\{0\}\cup F_{\mu_\Phi}(\R)\cup\{1\}.\ee
 
 In the case of \Lvy matrices, it has been proved in \cite[Theorem 1.3]{BDG} that $\mu_{\Phi}$ has no atoms (because it is absolutely continuous), so that  $E_\Phi=[0,1]$.

We introduce the eigenvalues $\lam_1\le \cd\le \lam_n$ of $A$ and  an orthogonal    matrix $U=[u_{ij}]$ \st 
$A=U\diag(\lam_1, \ld, \lam_n)U^*$. We  assume $U$ defined  in such a way that   the rows of the matrix $[|u_{ij}|]$ are exchangeable (this is possible\footnote{Such a matrix $U$ can be defined, for example, by choosing some orthogonal bases of all eigenspaces of $A$ with uniform distributions, independently with each other and independently of $A$ (given its eigenspaces of course).} because $A$ is invariant, in law, by conjugation by any permutation matrix). 
Then define the bivariate processes $$B^n_{s,t}:=\ff{\sqrt{n}}\sum_{\substack{1\le i\le ns\\ 1\le j\le nt}}(|u_{ij}|^2 -\ff{n})\qquad\trm{($0\le s \le 1$, \quad $0\le t \le 1$)}$$ and $$C^n_{s,\lam}:=\ff{\sqrt{n}}\sum_{\substack{1\le i\le ns\\ 1\le j\le n\ste \lam_j\le \lam}}(|u_{ij}|^2 -\ff{n})\qquad\trm{($0\le s \le 1$, \quad $\lam\in \R$)}.$$

The following theorem is the main result of this article. We endow $D([0,1]^2)$ and $D([0,1]\ti \R)$ with the  Skorokhod topology and $D([0,1]\ti E_\Phi)$ with the topology induced by the Skorokhod topology on $D([0,1]^2)$ by the projection map from $D([0,1]^2)$ onto $D([0,1]\ti E_\Phi)$ (see Section 4.1 of \cite{FloEigenvectors} for the  corresponding definitions).

%Recall that $E_\Phi\subset[0,1]$ has been defined at \eqre{defE_Phi} by the formula $E_\Phi:=\{0\}\cup F_{\mu_\Phi}(\R)\cup\{1\}$.
\beg{Th}\la{mainresult30513}As $n\to\infty$, the joint distribution of the processes $$(B^n_{s,t})_{(s,t)\in [0,1]\ti E_\Phi} \qquad \trm{ and }\qquad(C^n_{s,\lam})_{(s,\lam)\in [0,1]\ti \R}  $$ converges weakly to the joint  distribution of some centered Gaussian processes $$(B_{s,t})_{(s,t)\in [0,1]\ti E_\Phi} \qquad \trm{ and }\qquad(C_{s,\lam})_{(s,\lam)\in [0,1]\ti \R}  $$ vanishing on the boundaries of their domains and satisfying the relation \be\la{rel31513}B_{s, F_{\mu_\Phi}(\lam)}=C_{s,\lam}\ee for all $s\in [0,1]$, $\lam\in \R$.  Moreover, the process  $(B_{s,t})_{(s,t)\in [0,1]\ti E_\Phi}$ is continuous.
\en{Th}

\beg{rmk}Note that the limit of $B^n_{s,t}$ is only given here when $t\in E_\Phi$, \ie when $t$ is not in the ``holes" of $F_{\mu_\Phi}(\R)$. But as these holes result from the existence of some atoms in the limit spectral distribution of $A$, the variations  of $B^n_{s,t}$ when $t$ varies in one of these holes  may especially depend on the way we choose the columns of $A$ for eigenvalues with multiplicity larger than one.  By the results of \cite{cat-alain10}, in the case where the atoms of $\mu_\Phi$ result in atoms (with asymptotically same weight) of $\mu_n$,  the choice we made here should lead to a limit process $(B_{s,t})_{(s,t)\in [0,1]^2}$ which would interpolate  $(B_{s,t})_{(s,t)\in [0,1]\ti E_\Phi}$  with some    Brownian bridges in these ``holes", namely for when  $t\in  [0,1]\bck E_\Phi$.\en{rmk}

The following proposition     insures  that the $\ds\ff{\sqrt{n}}$ scaling in the definitions of $B^n_{s,t}$ and $C^n_{s,\lam}$ is the right one.

\beg{propo}\la{NonZeroCovariance}If the function $\Phi(z)$ of \eqre{2071216h33} is not linear in $z$, then    for any fixed $s\in (0,1)$, the covariance of the process $(B_{s,t})_{t\in E_\Phi}$   (hence also that of $(C_{s,\lam})_{\lam\in \R}$) is not identically null.
\en{propo}

\beg{rmk}One could wonder if the covariance might vanish uniformly on some compact in the $t$ variable, hence giving some support to the belief that the eigenvectors could behave more alike 
the eigenvectors of GUE for ``small'' eigenvalues (in the latter case the covariance should vanish).  Unfortunately, it does not seem that the covariance should be so closely related with the localization/delocalization properties of the eigenvectors (see   Remark \re{2441416h} instead). Indeed, 
let us consider  \Lvy matrices with $\alpha\in (1,2)$.   Their eigenvectors  are  delocalized \cite{charles_alice}, so that  one could expect the covariance of the process $(B_{s,t})_{t\in E_\Phi}$ to vanish.  This
is in contradiction with the fact that
 such matrices enter our model, hence have eigenvectors satisfying Theorem \re{mainresult30513} and Proposition \re{NonZeroCovariance}. \en{rmk}

%\beg{propo}\la{NonZeroCovarianceloc}If the function $\Phi(z)$ of \eqre{2071216h33} is not identically equal to $z$, then    for any  $s\in (0,1)$ and any $\lam<\lam'$,   $\var(C_{s,\lam'}-C_{s,\lam})\ne 0$.
%\en{propo}
%
%The following corollary follows directly from the previous proposition and \eqre{rel31513}.
%\beg{cor}If the function $\Phi(z)$ of \eqre{2071216h33} is not identically equal to $z$, then    for any  $s\in (0,1)$ and  any $0\le t<t'\le 1$ \st $\mu_{\Phi}((t,t'])>0$, $\var(B_{s,t'}-B_{s,t})\ne 0$.
%\end{cor}

To prove Theorem \re{mainresult30513}, a key step will be to prove the following proposition, which also allows to make the variance of the limiting processes in Theorem \re{mainresult30513} more explicit. 

Let us   define, for $z\in \C\bck \R$ and $s\in [0,1]$, \be\la{defXn(g,z)}X^n(s,z):=\ff{\sqrt{n}} \lf( \Tr (P_s \ff{z-A})- s_n \Tr \ff{z-A}\ri),\ee
where  $P_s$ denotes the diagonal matrix with diagonal entries $\one_{i\le ns}$ ($1\le i\le n$) and \be\la{111421h12}s_n:=\ff{n}\Tr P_s=\f{\lfl ns\rfl}{n}.\ee

 \begin{propo}\la{PropoCVST}
The distribution of the random process $$(X^n(s,z))_{s\in [0,1], z\in \C\bck\R}$$ converges weakly in the sense of finite marginals towards the distribution of a centered Gaussian process \be\la{3051317h15}(H_{s,z})_{s\in [0,1], z\in \C\bck\R}\ee with a covariance given by \eqre{591300h29}.
\end{propo} 

As it will appear from the proofs that the process $(C_{s,\lam})$ of Theorem \re{mainresult30513} and the process $(H_{s,z})$ from the previous proposition are linked by the formula \be\la{891312h22}  \int_{\lam\in \R}\f{ C_{s, \lam}}{(z-\lam)^2}\ud \lam  =-H_{s,z}\qquad\trm{($s\in [0,1]$, $z\in \C\bck\R$)},\ee
the covariance of $C_{s,\lambda}$ (hence of $B_{s,t}$ by \eqre{rel31513}) can be deduced from that of the process $H_{s,z}$ as follows (the proof of this proposition is a direct application of \eqre{891312h22} and of Formula \eqre{891312h18} of the Appendix).

% Let us define the slight modification $\tilde{C}_{s,\lam}$ of  $C_{s,\lam}$ by $$\tilde{C}_{s,\lambda}:=\f{1}{2}\bigg(\lim_{\tilde{\lam}\underset{<}{\to\lam}}C_{s,\tilde{\lam}}+C_{s, \lam}\bigg)$$
%(the difference between $C_{s,\lam}$ and $\tilde{C}_{s,\lam}$ appears only at the jumps in the second variable :   the discontinuities $\tilde{C}_{s,\lam}$ hare symmetrical discontinuities, whereas the function $C_{s, \lam}$ is \cadlag in its second variable). 

\beg{propo}\la{791316h37} For any $s,s'\in [0,1]$ and any $\lam, \lam'\in \R$ which are not atoms of $\mu_\Phi$, we have  \be\la{891312h29} \E[ {C}_{s,\lambda} {C}_{s',\lambda'}]=\ff{\pi^2}\lim_{\eta\downarrow 0}
\int_{-\infty}^\lambda\int_{-\infty}^{\lambda'}\mathbb E[ \Im \left( H_{s,E+i\eta}\right) \Im \left( H_{s',E'+i\eta}\right)]\ud E\ud E'\,.\ee When $\lam$ or $\lam'$ is an atom of $\mu_\Phi$, the covariance can be obtained using \eqre{891312h29} and the right continuity of $C_{s,\lam}$ in $\lam$.
\en{propo}

\section{Proof of Theorem \re{mainresult30513}}\la{secProof1}
We   introduce the cumulative distribution function \be\la{305131bis}F_n(\lam):=\ff{n}|\{j\ste \lam_j\le \lam\}|\ee of the empirical spectral law $\mu_n$ defined at \eqre{3051ter}.
We shall use the following formula several times: for all $s\in [0,1]$ and $\lam\in \R$, \be\la{305131}C^n_{s,\lam}=B^n_{s,F_n(\lam)}  .
\ee  

We know, by Lemma \re{lemCtightness} of the appendix, that  the sequence $(\op{distribution}(B^n))_{n\ge 1}$ is tight and has all its accumulation points   supported by the set of continuous functions on $[0,1]^2$. As $F_n$ converges to $F_{\mu_\Phi}$ in the Skorokhod topology, it follows that the sequences $$\tB^n:=(B^n_{s,t})_{(s,t)\in [0,1]\ti E_\Phi}\qquad \trm{ and }\qquad(C^n_{s,\lam}=B^n_{s,F_n(\lam)})_{(s,\lam)\in [0,1]\ti \R} $$
are tight in their respective spaces. 
To prove the theorem, it suffices to prove that  the sequence   $(\op{distribution}(\tB^n,C^n))_{n\ge 1}$ has   only one accumulation point (which is Gaussian centered, vanishing on the boundaries, supported by continuous functions as far as the first component  is concerned and satisfying \eqre{rel31513}). So let $((\tB_{s,t})_{(s,t)\in [0,1]\ti E_\Phi} , (C_{s,\lam})_{(s,\lam)\in [0,1]\ti \R})$ be a pair of random processes having for distribution such an accumulation point. By \eqre{305131}, we have  $$\tB_{s, F_{\mu_\Phi}(\lam)}=C_{s,\lam}$$ for all $s\in [0,1]$, $\lam\in \R$. Hence it suffices to prove that the distribution of $C$ is totally prescribed and Gaussian centered. 

First, let us note that one can suppose that along the corresponding   subsequence, the distribution  of   $((B^n_{s,t})_{(s,t)\in [0,1]^2},( C^n_{s,\lam})_{(s,\lam)\in [0,1]\ti\R})$ converges weakly to the distribution of a pair  $(B,C)$ of processes \st $B$ is continuous and vanishing on the boundary of $[0,1]^2$. The difference with what was supposed above is that now, $t$ varies in $[0,1]$ and not only  in $E_\Phi$. Again, by \eqre{305131}, we have  \be\la{rel31513dcta}B_{s, F_{\mu_\Phi}(\lam)}=C_{s,\lam}\ee for all $s\in [0,1]$, $\lam\in \R$. Hence the process $C$ is continuous in $s$ and continuous in $\lam$ at any $\lam$ which is not an atom of the (non random) \pro measure  $\mu_\Phi$. 
Hence it follows from Lemma \re{lemcauchytransfo} of the appendix that it suffices to prove that the distribution of the process $$\lf(X(s,z):=\int_{\lam\in \R}\f{ C_{s,\lam}}{(z-\lam)^2}\ud \lam\ri)_{s\in [0,1], z\in \C\bck\R}$$ is totally prescribed (and Gaussian centered). This distribution is the limit distribution, along our subsequence, of the  process \be\la{16131h43}\lf(\int_{\lam\in \R}\f{ C^n_{s, \lam}}{(z-\lam)^2}\ud \lam\ri)_{s\in [0,1], z\in \C\bck\R}.\ee But 
   by  Lemma \re{lem3051317h13} below, the process of \eqre{16131h43} is simply  (the opposite of) the process $(X^n(s,z))_{s,z}$, defined above at 
\eqre{defXn(g,z)}. 
As Proposition 
\re{PropoCVST}
states that (regardless of the subsequence considered) the distribution of the process $(X^n(s,z))_{s,z}$ converges weakly to a Gaussian centered limit, this concludes the proof of Theorem \re{mainresult30513}.

\beg{lem}\la{lem3051317h13} For any $s\in [0,1]$ and any $z\in \C\bck\R$, we have \be\la{305132} \int_{\lam\in \R}\f{ C^n_{s, \lam}}{(z-\lam)^2}\ud \lam  =-X^n(s,z).\ee  
\en{lem}
 
 \bpr Let us introduce, for $s\in [0,1]$,    the  random signed measure $ \nu_{s,n} $ on $\R$ defined by  $$\nu_{s, n}:= \ff{\sqrt{n}}\sum_{1\le i\le ns}\sum_{j=1}^n \lf(|u_{ij}|^2-\ff{n}\ri) \del_{\lam_j}.$$ Then for any $\lam\in\R$, $\ds 
 C^n_{s,\lam}=\nu_{s,n}((-\infty, \lam]).$ 
Moreover, by Fubini's theorem, we know that for any finite signed measure $m$ on $\R$,   \be\la{67136h15}\int_{\lam\in \R}\f{m((-\infty,\lam])}{(z-\lam)^2}\ud \lam=-\int_{\lam\in \R}\f{\ud m(\lam)}{z-\lam}.\ee
Hence $$\int_{\lam\in \R}\f{ C^n_{s, \lam}}{(z-\lam)^2}\ud \lam =-\int_{\lam\in \R}\f{\ud\nu_{s,n}(\lam)}{z-\lam}.$$
On the other hand,  we have   \beq X^n(s,z)&=&\ff{\sqrt{n}} \lf(\sum_{1\le i\le ns}\lf(\ff{z-A}\ri)_{ii}-\ff{n}\sum_{1\le i\le ns}\sum_{j=1}^n\ff{z-\lam_j}\ri)\\
&=&\ff{\sqrt{n}} \lf(\sum_{1\le i\le ns}\sum_{j=1}^n |u_{ij}|^2\ff{z-\lam_j}-\ff{n}\sum_{1\le i\le ns}\sum_{j=1}^n\ff{z-\lam_j}
 \ri)\\ 
 &=&\ff{\sqrt{n}}  \sum_{1\le i\le ns}\sum_{j=1}^n \lf(|u_{ij}|^2-\ff{n}\ri)\ff{z-\lam_j} \\
 &=& \int_{\lam\in \R}\f{\ud\nu_{s,n}(\lam)}{z-\lam}.
 \eeq This concludes the proof. 
 \epr

\section{Proof of Proposition 
\re{PropoCVST}}\la{sec:PropoCVST}

To prove Proposition 
\re{PropoCVST}, one needs to prove that the distribution of any linear combination of the $X^n(s,z)$'s ($s\in [0,1]$, $z\in \R$) converges weakly. For $s=0$ or $1$, $\nu_{s,n}$ is null, as $X^n(s,z)$, hence we can focus on $s\in (0,1)$. Any such linear combination can be written $$M^n:=\sum_{i=1}^p\al_i  X^n({s_i},z_i),$$ for some $\al_i$'s in $\C$, some $s_i$'s in $[0,1]$ and some complex non real numbers $z_i$.

We want to prove that  $M^n$ converges in law to a certain complex centered Gaussian variable.  
We are going  to use the CLT for martingale differences stated at  Theorem \re{thconvmart} of the appendix. Indeed,  for $\F_k^{\, n}$ the $\si$-algebra generated by the first $k\ti k$ upper-left corner of the symmetric matrix $A$, the sequence  $(M_k^n:=\E[M^n|\F_k^{\, n}])_{k=0,\ld, n}$ is a centered martingale  (to see that it is centered, just use the fact that as $A$ is invariant, in law, by conjugation by any permutation matrix, for all $z$, the expectation of $(\ff{z-A})_{jj}$ does not depend on $j$).

Then,   denoting $\E[\,\cdot\,|\mc{F}_k^{\, n}]$ by $\E_k$, and defining $$Y_k:=(\E_k-\E_{k-1})(M^n)$$ (which depends implicitely on $n$),  we need to prove 
that for any $\eps>0$, 
\be\la{16132h34}L^n(\eps):=\sum_{k=1}^n \E(|Y_k|^2\one_{|Y_k|\ge \eps})\ninf 0,\ee
and that the sequences $$\sum_{k=1}^n \E_{k-1}(|Y_k|^2)\qquad\trm{ and }\qquad \sum_{k=1}^n \E_{k-1}(Y_k^2)$$   converge  in \pro
 towards some deterministic limits.
 As $\ovl{X^n(s,z)}=X^n(s,\ovl{z})$, it is in fact enough to fix $s,s'\in (0,1)$ and $z,z'\in \C\bck\R$ and to  prove that for 
 \be\la{defY_kfinal} Y_k:=(\E_k-\E_{k-1})(X^n(s,z))\qquad\trm{ and }\qquad Y_k':=(\E_k-\E_{k-1})(X^n(s',z')),
 \ee 
 we have \eqre{16132h34} for any $\eps>0$ and that 
  $$\sum_{k=1}^n \E_{k-1}(Y_kY_k')$$ converges in \pro towards a deterministic constant.
We introduce the notation $$G:=\ff{z-A}\qquad \trm{ and }\qquad G':=\ff{z'-A}.$$ 

Recall that  $P_s$ denotes the diagonal matrix with diagonal entries $\one_{i\le ns}$ ($1\le i\le n$).
Let $A^{(k)}$ be the symmetric matrix with size $n-1$ obtained by removing the $k$-th row and the $k$-th column of $A$.  The matrix $P_s^{(k)}$ is defined in the same way out of $P_s$. Set $G^{(k)}:=\ff{z-A^{(k)}}$. 
Note that $\E_k G^{(k)}=\E_{k-1}G^{(k)}$, so that $Y_k$, which is equal to $\ds \ff{\sqrt{n}} (\E_k-\E_{k-1})\lf( \Tr (P_sG)- s_n\Tr G\ri)$, can be rewritten  \beq Y_k
&=&\ff{\sqrt{n}} (\E_k-\E_{k-1})\lf( (\Tr (P_sG)-\Tr (P_s^{(k)}G^{(k)}))-s_n (\Tr G-\Tr G^{(k)})\ri)
\eeq

Then, \eqre{16132h34} is obvious by Formula \eqre{161312h} of the appendix (indeed, $L^n(\eps)$ is null for $n$ large enough).
Let us now  apply  Formula  \eqre{104138h} of the appendix. We get
\be\la{161312h36}Y_k=\ff{\sqrt{n}} (\E_k-\E_{k-1})\lf(\f{\one_{k\le ns}-s_n+\ba_k^*G^{(k)}(P_s^{(k)}-s_n)G^{(k)}\ba_k}{z-a_{kk}  -\ba_k^* G^{(k)}\ba_k}
\ri).
\ee

    Following step by step Paragraph 3.2 of  \cite{ACFTCL}, one can neglect the non diagonal terms in the expansions of the quadratic forms in \eqre{161312h36}, \ie replace $Y_k$ by $\ds \ff{\sqrt{n}}(\E_k-\E_{k-1})\lf(f_k \ri)$, with \be\la{deff_k61313h17}f_k :=f_k(z,s)= \f{\one_{k\le ns}-s_n+\sum_j\ba_k(j)^2\{G^{(k)}(P_s^{(k)}-s_n)G^{(k)}\}_{jj} }{z  -\sum_j\ba_k(j)^2 G^{(k)}_{jj}}
 . \ee
%    The proof of this statement  mimics the one of   Section 3.2 of \cite{ACFTCL}: it is enough to show that
%    $$\sum_{i\neq \ell}a_{ki}\bar a_{k\ell} (G^{(k)}(P_s^{(k)}-s_n)G^{(k)})_{i\ell} $$
%    is  small in probability: this is a direct consequence of the fact that  $n\mathbb \E[a_{k\ell}^2]$ is bounded or that $ a_{k\ell}$ is in the domain of attraction of an $\alpha$-stable law, together with the
%    uniform boundedness of the spectrums of $(G^{(k)}(P_s-s_N)G^{(k)})$ and $(G^{(k)})^2$ and of  $f_k $.
In other words,  \be\la{etavril2013}\sum_{k=1}^n \E_{k-1}(Y_kY_k')= \ff{n}\sum_{k=1}^n\E_{k-1}\lf[(\E_k-\E_{k-1})\lf(f_k \ri)(\E_k-\E_{k-1})\lf(f_k' \ri)\ri] 
    +o(1),\ee where $f_k'$ is defined as $f_k$ in \eqre{deff_k61313h17}, replacing the function $s$ by $s'$ and $z$ by $z'$.

    Let us denote by $\E_{\ba_k}$ the expectation with respect to the randomness of the $k$-th column of $A$ (\ie the conditional expectation with respect to the $\si$-algebra generated by the $a_{ij}$'s \st $k\notin \{i,j\}$). Note that $\E_{k-1}=\E_{\ba_k}\circ\E_k=\E_{k}\circ \E_{\ba_k}$,  hence 
  \be\la{38128h50vect} \E_{k-1}\lf[(\E_k-\E_{k-1})\lf(f_k \ri)(\E_k-\E_{k-1})\lf(f_k' \ri)\ri] =\E_{k}[\E_{\ba_k}(f_k\ti f_k'')]-\E_{k}\E_{\ba_k}f_k\ti \E_{k}\E_{\ba_k}f_k',\ee
    where $f_k''$ is defined as $f_k'$  replacing the matrix $A$ by the matrix   
 \be\la{88123}
		  A'=[ a_{ij}']_{1\le i,j\le N}
	\ee
  defined by  the the fact that the   $ a_{ij}'$'s \st $i>k$ or  $j>k$ are i.i.d. copies of $a_{11}$ (modulo the fact that $A'$ is symmetric),     independent of $A$ and 
 for all other pairs $(i,j)$, 
$ a_{ij}'=a_{ij}$.

  For each $s\in (0,1)$ let us define $\C^2_s$ to be  the set of pairs $(z,\tz)$ of complex numbers \st $$(\Im z> 0\trm{ and }-\f{\Im z}{1-s}<\Im \tz<\f{\Im z}{s})\qquad\trm{ or }\qquad(\Im z< 0\trm{ and }\f{\Im z}{s}<\Im \tz<-\f{\Im z}{1-s}).$$
Note that $\C^2_s$ is the set of pairs $(z,\tz)$ of complex numbers \st $\Im z\ne 0$ and both $\Im(z+(1-s)\tz)$ and $\Im(z-s\tz)$ have the same sign as $\Im z$. At last, in the next lemma, $$\pa_{\tz}=\f{\pa}{\pa \tz}$$ is not to be taken for the usual notation $\pa_{\bar{z}}$.
    
    \begin{lem}\la{204131} For any fixed $z\in \C\bck\R$ and any fixed $s\in (0,1)$,  as $n,k\lto\infty$ in such a way that $k/n\lto u\in [0,1]$, we have the convergence in probability  $$\lim_{N\ra\infty} \E_{\ba_k}[f_k(z,s)]=L_u(z,s):= -\int_0^{+\infty}\ff{t}\pa_{\tz, \tz=0}e^{i\sgn_z t(z+\tz(\one_{u\le s}-s))}e^{\rho_{z,\tz, s}(t)}\ud t,$$   where for $s\in (0,1)$ fixed,  $(z,\tz,t)\longmapsto\rho_{z,\tz, s}(t)$ is the unique  function defined on   $\C^2_s\ti\R_+$, analytic in its two first variables and continuous in its third one, taking values  into $\{z\in \C\ste \Re z\le 0\}$,  solution of 
$$\rho_{z,\tz, s}(t)=t\int_0^\infty g(ty) (s e^{iy\sgn_z  \tz}+(1-s)) e^{iy\sgn_z(z-s \tz)}e^{ \rho_{z,\tz, s}(y)} \ud y$$
where $g$ is the function introduced at \eqre{hypcalcconv}.\end{lem}

    \bpr
    We use the fact that   for $z \in \C \bck\R$,  	\be\la{rep}
		\ff{z}=-i\op{sgn}_z\ti  \int_{0}^{+\infty}e^{\op{sgn}_z itz}\ud t,
	\ee 
	 where  $\op{sgn}_z$ has been defined above by $\op{sgn}_z=\op{sgn}(\Im z)$. 
	 Hence by \eqre{deff_k61313h17},   $$f_k=-i\sgn_z\int_0^{+\infty}\{\one_{k\le ns}-s_n+\sum_j(G^{(k)}(P_s^{(k)}-s_n)G^{(k)})_{jj}\ba_k(j)^2\}e^{i\sgn_z t(z  -\sum_j(G^{(k)})_{jj}\ba_k(j)^2)}\ud t.
$$
Let us define   $\ds G^{(k)}(\cdot,\cdot)$   on $\C_{s_n}^2$ by \be\la{111421h01}\ds G^{(k)}(z,\tz):=\ff{z+\tz(P_s^{(k)}-s_n)-A^{(k)}}\ee (note that $\ds G^{(k)}(\cdot,\cdot)$ is well defined by the remark following the definition of $\C^2_s$). Then   for any fixed $z\in \C\bck\R$, 
$$G_k(z)(P_s^{(k)}-s_n)G_k(z)= -\pa_{\tz, \tz=0}G^{(k)}(z,\tz).$$ 
Hence 
$$\{\one_{k\le ns}-s_n+\sum_j(G^{(k)}(P_s^{(k)}-s_n)G^{(k)})_{jj}\ba_k(j)^2\}e^{i\sgn_z t(z  -\sum_j(G^{(k)})_{jj}\ba_k(j)^2)}$$ $$=\ff{it\sgn_z}\pa_{\tz, \tz=0}e^{i\sgn_z t\{z+\tz(\one_{k\le ns}-s_n)-\sum_j(G^{(k)}(z,\tz)_{jj}\ba_k(j)^2\}}$$
and 
\be\la{161317h36}f_k=-\int_0^{+\infty}\ff{t}\pa_{\tz, \tz=0}e^{i\sgn_z t\{z+\tz(\one_{k\le ns}-s_n)-\sum_jG^{(k)}(z,\tz)_{jj}\ba_k(j)^2\}}\ud t.
\ee
%(as $z\in \C\bck \R$ is fixed and for each $j$, $-G^{(k)}(z,\tz)_{jj}$ has imaginary part with the same sign as $z$ for $\tilde{z}$ small enough, the differentiation under the integral is justified).
Let us now compute $\E_{\ba_k}(f_k)$. One can  permute $\E_{\ba_k}$ and $\int_0^{+\infty}$ because  $z\in \C\bck \R$ is fixed and for each $j$, $-G^{(k)}(z,\tz)_{jj}$ has imaginary part with the same sign as $z$ for $\tilde{z}$ small enough. Hence  for $\phi_n$  defined as in \eqref{2071216h23}  by $\phi_n(\lam)=\E e^{-i\lam a_{11}^2}$, we have 
\beq \E_{\ba_k}(f_k)
&=&- \int_0^{+\infty}\ff{t}\pa_{\tz, \tz=0}e^{i\sgn_z t(z+\tz(\one_{k\le ns}-s_n)}\prod_j \phi_n(\sgn_z tG^{(k)}(z,\tz)_{jj}) \ud t\eeq  Now, by \eqre{2071216h33}, we have the uniform convergence on compact sets $n(\phi_n-1)\lto \Phi$ as $n\to\infty$.
As $\Re(i\sgn_z z)<0$, the integrals are well dominated at infinity. Moreover,      the integral $$\int_0^{+\infty}\ff{t}\pa_{\tz, \tz=0}e^{i\sgn_z t(z+\tz(\one_{k\le ns}-s_n))}e^{\ff{n}\sum_j\Phi(\sgn_ztG^{(k)}(z,\tz)_{jj})}\ud t$$ is well converging at 
 the origin as the derivative in $\tilde z$ is of order $t$.  Indeed,  $G^{(k)} (z,\tilde z)_{jj}$   takes its values in $\C^-$ and is  uniformly bounded, 
 and  $\Phi$ is analytic on   $\mathbb C^-$. By Lemma \re{2541400h15}, it follows that    
\beq \E_{\ba_k}(f_k)
&=&- \int_0^{+\infty}\ff{t}\pa_{\tz, \tz=0}e^{i\sgn_z t(z+\tz(\one_{k\le ns}-s_n))}e^{\ff{n}\sum_j\Phi(\sgn_ztG^{(k)}(z,\tz)_{jj})}\ud t +o(1)\,.
\eeq   
We therefore basically need to compute the asymptotics of $$\rho_{z,\tz, s}^n(t):=\ff{n}\sum_{j}\Phi(\sgn_ztG^{(k)}(z,\tz)_{jj}).$$
Note that by definition of $\Phi$, for any $\lam\in \ovl{\C^-}$, $\Re\Phi(\lam)\le 0.$
Thus $\rho_{z,\tz, s}^n(t)$ is analytic in $z\in \mathbb C\backslash \mathbb R$, and uniformly bounded on   compact subsets of $\C\backslash \mathbb R$ and takes values in $\{z\in \C\ste \Re z \le 0\}$. By Montel's theorem, all limit points of this function for   uniform convergence on compact subsets will satisfy the same property. 
Now, notice   by Schur complement formula and the removal of the non diagonal terms (Lemma 7.7 of \cite{ACFTCL} again), that for $n\gg 1$, 
\be\la{244144h1}G^{(k)}(z,\tz)_{jj}= \frac{1}{ z+\tz(\one_{k\le ns} -s_n)-\sum_\ell a_{j\ell}^2 G^{(k,j)}(z,\tz)_{\ell\ell}}+o(1)\ee
where $G^{(k,j)}$ is the resolvent where two rows and columns have been suppressed.
We can now proceed  to write that
by invariance of the law of $A$ by conjugation by permutation matrices,  for all $j$,  $$\E[\Phi(\sgn_ztG^{(k)}(z,\tz)_{jj})]=\beg{cases}\E[\Phi(\sgn_ztG^{(k)}(z,\tz)_{11})]&\trm{ if $j\le ns$,}\\ \\ 
 \E[\Phi(\sgn_ztG^{(k)}(z,\tz)_{nn})]&\trm{ if $j> ns$,}\en{cases}$$ 
so that by concentration arguments, see \cite[Appendix]{ACFTCL}, $\rho_{z,\tz, s}^n(t)$ self-averages and  for $n\gg 1$, with very large probability, 
\beq \rho_{z,\tz, s}^n(t)&=& \E[\frac{1}{n}\sum_{j} \Phi(\sgn_ztG^{(k)}(z,\tz)_{jj})]+o(1)\\
&=& s_n \E[\Phi(\sgn_z tG^{(k)}(z,\tz)_{11})]+(1-s_n) \E[\Phi(\sgn_z tG^{(k)}(z,\tz)_{nn})]+o(1).\eeq
On the other side, using \eqre{244144h1},   the function $g$ introduced in the hypothesis at \eqre{hypcalcconv} and a change of variable $y\to y/t$, we have (using Lemma \re{2541400h15} twice)
\begin{eqnarray}\nonumber
\E[\Phi(\sgn_ztG^{(k)}(z,\tz)_{11})]&=&t\int_0^\infty g(ty) e^{iy \sgn_z (z+\tz(1-s_n))} \prod_j\phi_n( y \sgn_zG^{(k,1)}(z,\tz)_{jj})\ud y\\ \la{lem204131eq1}
&=& t\int_0^\infty g(ty) \exp(iy \sgn_z (z+\tz(1-s_n))) e^{ \rho_{z,\tz, s}^n(y)} \ud y+o(1)\\
\nonumber\E[\Phi(\sgn_z tG^{(k)}(z,\tz)_{nn})]&=&t\int_0^\infty g(ty) e^{iy\sgn_z( z-s_n \tz)} \prod_j \phi_n( \sgn_z y G^{(k,1)}(z,\tz)_{jj})\ud y\\
 \la{lem204131eq2}&=& t\int_0^\infty g(ty) \exp(iy\sgn_z( z-s_n \tz)) e^{ \rho_{z,\tz, s}^n(y)} \ud y
+o(1)\end{eqnarray}
so that we deduce that
the limit points $\rho_{z,\tz, s}(t)$ of $\rho_{z,\tz, s}^n(t)$ satisfy
$$\rho_{z,\tz, s}(t)=t\int_0^\infty g(ty) (s e^{iy\sgn_z  \tz}+(1-s)) e^{iy\sgn_z(z-s \tz)}e^{ \rho_{z,\tz, s}(y)} \ud y.$$

Let us now prove that for each fixed $s\in (0,1)$, 
  there exists a unique function satisfying this equation and the conditions stated in the lemma. 
So let us suppose that we have two solutions $\rho_{z,\tz, s}(t)$ and $\tilde{\rho}_{z,\tz, s}(t)$ with non positive real parts. Then $$\Del_{z,\tz}(t):=\rho_{z,\tz, s}(t)-\tilde{\rho}_{z,\tz, s}(t)$$ satisfies $$\Del_{z,\tz}(t) =t\int_0^\infty g(ty) (s e^{iy\sgn_z  \tz}+(1-s)) e^{iy\sgn_z(z-s \tz)}(e^{ \rho_{z,\tz, s}(y)} -e^{ \tilde{\rho}_{z,\tz, s}(y)})\ud y.$$ Let  $\del(z,\tz):=\min\{\sgn_z \Im(z+(1-s)\tz),\;\sgn_z\Im(z-s\tz) \}>0.$  We have $$ |(s e^{iy\sgn_z  \tz}+(1-s)) e^{iy\sgn_z(z-s \tz)}|\le e^{-\del(z,\tz) y},$$ hence 
  \beq |\Del_{z,\tz}(t)| &\le &t\int_0^\infty |g(ty)| e^{-\del(z,\tz) y} |\Del_{z,\tz}(y )|\ud y\eeq
Thus by the hypothesis made on $g$ at \eqre{ConditionOng}, 
 \beq |\Del_{z,\tz}(t)| &\le &Kt^{\ga+1}\underbrace{\int_0^\infty y^\ga e^{-\del(z,\tz) y} |\Del_{z,\tz}(y )|\ud y}_{:=I_1(z,\tz)}+Kt^{\ka+1}\underbrace{\int_0^\infty y^\ka e^{-\del(z,\tz) y} |\Del_{z,\tz}(y )|\ud y}_{:=I_2(z,\tz)}\eeq
 It follows that the numbers $I_1(z,\tz)$ and $I_2(z,\tz)$ defined above satisfy \begin{eqnarray*}
 I_1(z,\tz) &\le&  K\bigg( I_1(z,\tz) \int_{0}^\infty y^{2\ga+1}e^{-\del(z,\tz) y}\ud y + I_2(z,\tz) \int_{0}^\infty y^{\ga+ \ka+1}e^{-\del(z,\tz) y}\ud y  \bigg),\\
  I_2(z,\tz) & \le &  K\bigg( I_1(z,\tz) \int_{0}^\infty y^{\ga+ \ka+1}e^{-\del(z,\tz) y}\ud y + I_2(z,\tz) \int_{0}^\infty y^{2\ka+1}e^{-\del(z,\tz) y}\ud y  \bigg).\end{eqnarray*}
   For $\del(z,\tz)$ large enough, 
   the integrals above are all strictly less that $\frac 1 {4K}$,
    so  
   $I_1(z,\tz)=I_2(z,\tz)=0$.
It follows that for $\Im z$ large enough and $\Im \tz$ small enough, both solutions   coincide. By analytic continuation, unicity follows.
\epr

        Getting back to \eqre{etavril2013} and \eqre{38128h50vect}, we shall now,  as in \cite{ACFTCL},  analyze
         \be\la{11142}L_k^n(s,z;s',z'):=\E_{\ba_k}(f_k\ti f_k'').\ee
Let us first  define the measure  \be\la{2612131}\tilde\tau:=\tau+\delta_0\otimes\mu+\mu\otimes \delta_0\ee on $(\R^+)^2$     for $\tau$ and $\mu$ the measures introduced at   \eqre{phias} or at Remark \re{22101318h17}.  
We always have, for $x,y\in \C^+$,
\be\label{phias221013}
\Phi(x^{-1}+y^{-1})=\iint_{(\R_+)^2}e^{i(xv+yv')} \ud\tilde{\tau}(v,v')\ee
         \begin{lem}\label{defL2} Let us fix $s_1,s_2\in (0,1)$.  As $k,n\lto\infty$ in such a way that  $k/n$ tends to $u\in [0,1]$, the quantity $L_k^n(s_1,z;s_2,z')$ defined at \eqre{11142} converges in \pro  to the deterministic limit 
 $$L_u(s_1,z;s_2,z'):=$$ $$ \iint_{\R_+^2} \partial_{\tz,\tz=0}\partial_{\tz',\tz'=0} e^{i\sgn_z t(z+\tz(\one_{u\le s_1}-s_1))+i\sgn_{z'} t'(z'+\tz'(\one_{u\le s_2}-s_2))+\rho_u(s_1,t,z,\tz;s_2,t',z',\tz')} 
       \f{\ud t\ud t'}{tt'}
$$ where the function $$(t,z,\tz,t',z',\tz')\in \R_+\ti\C^2_{s_1}\ti\R_+\ti \C^2_{s_2} \longmapsto \rho_u(s_1,t,z,\tz;s_2,t',z',\tz')$$ is characterized as follows : 
    \be\la{2441417h}\rho_u(s_1,t,z,\tz;s_2,t',z',\tz')=\rho_u(s_2,t',z',\tz';s_1,t,z,\tz)\ee and 
      if, for example,  $s_1\le s_2$,  then for $\gamma_1=s_1,\gamma_2=s_2-s_1$, $\gamma_3=1-s_2$ and $\tilde\tau=\tau+\delta_0\otimes\mu+\mu\otimes \delta_0$,
          \beq \rho_u(s_1,t_1,z_1,\tz_1;s_2,t_2,z_2,\tz_2)&=&u\sum_{\bet=1}^3\ga_\bet \iint_{\R_+^2}  e^{ \sum_{r=1,2}
\sgn_{z_r}\frac{iv_r}{t_r} \{z_r+\tz_r(\one_{\bet\le r}-s_r)\}}\ti \\ \nonumber&& \qquad e^{
 \rho_u(s_1,\f{v_1}{t_1},z_1,\tz;s_2,\f{v_2}{t_2},z_2,\tz_2)}\ud\tilde{\tau}(v,v')+\\ &&
 \sum_{r=1,2}t_r\int_0^\infty g(t_ry) \{(s_r-u)^+ e^{iy\sgn_{z_r}  \tz_r}+1-\max(s_r,u)\} 
 \\
 \nonumber&&\qquad e^{iy\sgn_{z_r}(z_r-s_r \tz_r)}e^{ \rho_{z_r,\tz_r, s_r}(y)} \ud y
 \eeq
 (the characterization of $\rho_u(s_1,t,z,\tz;s_2,t',z',\tz')$ when $s_2\le s_1$ can be deduced from the previous equation and \eqre{2441417h}).
 \end{lem}
 
 \bpr Of course, $L_k^n(s_1,z;s_2,z')=L_k^n(s_2,z';s_1,z)$. Let us suppose for example that $s_1\le s_2$.  We use the definition of $G^{(k)}(z,\tz)$ given at \eqre{111421h01} for $s$ replaced by $s_1$ and define in the same way,  for  $(z',\tz)\in \C_{s_{2,n}}^2$, $$ {G'}^{(k)}(z',\tz):= \ff{z'+\tz(P_{s_2}^{(k)}-s_{2,n})-{A'}^{(k)}}$$ with $s_{i,n}:= \f{\lfl ns_i\rfl}{n}$ ($i=1,2$).

       First, recall the following formula for $f_k
       $ established at Equation \eqre{161317h36}:  
        $$f_k=-\int_0^{+\infty}\ff{t}\pa_{\tz, \tz=0}e^{i\sgn_z t\{z+\tz(\one_{k\le ns}-s_n)-\sum_{j\ne k}G^{(k)}(z,\tz)_{jj}\ba_k(j)^2\}}\ud t.$$ 
       In the same way, we find 
        $$f_k''=-\int_0^{+\infty}\ff{t'}\pa_{\tz, \tz=0}e^{i\sgn_{z'} t'\{z'+\tz(\one_{k\le ns'}-s_n')-\sum_{j\ne k}{G'}^{(k)}(z',\tz)_{jj}\ba_k'(j)^2\}}\ud t.$$ 
        As the $\ba_k(j)$ and the $\ba_k'(j)$ are identical when $j\le k$ and independent when $j>k$,  we have 
        \begin{eqnarray*}  L_k^n(s_1,z_1;s_2,z_2) &=&
       \iint_{\R_+^2}        \partial_{\tz_1,\tz_1=0}\partial_{\tz_2,\tz_2=0}   e^{i\sgn_{z_1} t_1\{z_1+\tz_1(\one_{k\le ns_1}-s_{1,n})\}+i\sgn_{z_2} t_2\{z_2+\tz_2(\one_{k\le ns_2}-s_{2,n})\}}\ti\\
        &&\prod_{j\le k}\phi_n( \sgn_{z_1} t_1 G^{(k)}(z_1,\tz_1)_{jj}+\sgn_{z_2}t_2  {G'}^{(k)}(z_2,\tz_2)_{jj})\ti\\
        &&\prod_{j> k}\phi_n( \sgn_{z_1} t_1 G^{(k)}(z_1,\tz_1)_{jj})\phi_n(\sgn_{z_2}t_2  {G'}^{(k)}(z_2,\tz_2)_{jj})\f{\ud t_1\ud t_2}{t_1t_2}
        \eeq Then using the usual off-diagonal terms removal and Lemma \re{2541400h15}, we get 
              \begin{eqnarray*}  L_k^n(s_1,z_1;s_2,z_2)  &=& \iint_{\R_+^2}    \partial_{\tz_1,\tz_1=0}\partial_{\tz_2,\tz_2=0} 
        e^{\sum_{r=1,2}i\sgn_{z_r} t_r\{z_r+\tz_r(\one_{k\le ns_r}-s_{r,n})\}} \ti\\ &&
        \exp(\rho_k^n(s_1,t_1,z_1,\tz_1;s_2,t_2,z_2,\tz_2))\f{\ud t_1\ud t_2}{t_1t_2}+o(1)
        \end{eqnarray*}
        with 
        \begin{eqnarray*}
           \rho_k^n(s_1,t_1,z_1,\tz_1;s_2,t_2,z_2,\tz_2)&:=&\frac{1}{n}\sum_{j\le  k} \Phi( \sgn_{z_1} t_1 G^{(k)}(z_1,\tz_1)_{jj}+\sgn_{z_2}t_2  {G'}^{(k)}(z_2,\tz_2)_{jj})\\ \la{21142}\\ &+&  \frac{1}{n}\sum_{j> k} \{\Phi( \sgn_{z_1} t_1 G^{(k)}(z_1,\tz_1)_{jj})+\Phi(\sgn_{z_2}t_2  {G'}^{(k)}(z_2,\tz_2)_{jj})\}
\end{eqnarray*}
By the Schur formula, \eqre{phias221013} and  the off-diagonal terms removal, we have \begin{eqnarray}\la{11141}&&
 \Phi( \sgn_{z_1} t_1 G^{(k)}(z_1,\tz_1)_{jj}+\sgn_{z_2}t_2  {G'}^{(k)}(z_2,\tz_2)_{jj})=\\ \nonumber\\  \nonumber
&& \iint_{\R_+^2}  e^{ 
\sgn_z\frac{iv_1}{t_1} \{z_1+\tz_1(\one_{j\le ns_1}-s_1)-\sum_{\ell\notin\{k,j\}} a_{j\ell}^2 G^{(k,j)}(z_1,\tz_1)_{\ell\ell}\}}  \\ \nonumber &&\ti e^{ 
\sgn_{z_2}\frac{iv_2}{t_2} \{z_2+\tz_2(\one_{j\le ns_2}-s_2)-\sum_{\ell\notin\{k,j\}} (a'_{j\ell})^2 {G'}^{(k,j)}(z_2,\tz_2)_{\ell\ell}\}} \ud\tilde{\tau}(v_1,v_2)+o(1).\eeqy
By  \eqre{11141} 
and using    the concentration around the expectation, we get, for $j<k$,  \beqy\la{312131}  && \Phi( \sgn_{z_1} t_1 G^{(k)}(z_1,\tz_1)_{jj}+\sgn_{z_2}t_2  {G'}^{(k)}(z_2,\tz_2)_{jj}) =o(1)+\\
\nonumber&&  \iint_{\R_+^2}  e^{ 
\sgn_{z_1}\frac{iv_1}{t_1} \{z_1+\tz_1(\one_{j\le ns_1}-s_1)\}+ 
\sgn_{z_2}\frac{iv_2}{t_2} \{z_2+\tz_2(\one_{j\le ns_2}-s_2)\}}%\ti \\ \nonumber&&
 e^{\rho_k^n(s_1,\f{v_1}{t_1},z_1,\tz;s_2,\f{v_2}{t_2},z_2,\tz_2)}\ud\tilde{\tau}(v,v').
 \end{eqnarray}
 Now, using   the proof of Lemma \re{204131} (especially \eqre{lem204131eq1} and \eqre{lem204131eq2}) and the fact that $k/n\lto u$, we get
 \beqy\nonumber\frac{1}{n}\sum_{j> k} \Phi( \sgn_{z_1} t_1 G^{(k)}(z_1,\tz_1)_{jj})&=&o(1)+(s_1-u)^+\E(\Phi(\sgn_{z_1} t_1 G^{(k)}(z_1,\tz_1)_{11})+\\ \nonumber&&(1-\max(s_1,u))\E(\Phi(\sgn_{z_1} t_1 G^{(k)}(z_1,\tz_1)_{nn})\\
 \la{111400h26}&=& o(1)+t_1\int_0^\infty g(t_1y) \{(s_1-u)^+ e^{iy\sgn_{z_1}  \tz_1}+1-\max(s_1,u)\} 
 \\
 \nonumber&&\qquad e^{iy\sgn_{z_1}(z_1-s_1 \tz_1)}e^{ \rho_{z_1,\tz_1, s_1}(y)} \ud y
 \eeqy
 In the same way,  \beq \frac{1}{n}\sum_{j> k} \Phi( \sgn_{z_2} t_2 G^{(k)}(z_2,\tz_2)_{jj})
 &=& o(1)+t_2\int_0^\infty g(t_2y) \{(s_2-u)^+ e^{iy\sgn_{z_2}  \tz_2}+1-\max(s_2,u)\} 
 \\
 \nonumber&&\qquad e^{iy\sgn_{z_2}(z_2-s_2 \tz_2)}e^{ \rho_{z_2,\tz_2, s_2}(y)} \ud y
 \eeq
 Summing up, we get that any  limit point  $\rho_u(s_1,t_1,z_1,\tz_1;s_2,t_2,z_2,\tz_2)$ of $$\rho_k^n(s_1,t_1,z_1,\tz_1;s_2,t_2,z_2,\tz_2)$$  satisfies 
 \beq \rho_u(s_1,t_1,z_1,\tz_1;s_2,t_2,z_2,\tz_2)&=&u\sum_{\bet=1}^3\ga_\bet \iint_{\R_+^2}  e^{ \sum_{r=1,2}
\sgn_{z_r}\frac{iv_r}{t_r} \{z_r+\tz_r(\one_{\bet\le r}-s_r)\}}\ti \\ \nonumber&& \qquad e^{
 \rho_u(s_1,\f{v_1}{t_1},z_1,\tz;s_2,\f{v_2}{t_2},z_2,\tz_2)}\ud\tilde{\tau}(v,v')+\\ &&
 \sum_{r=1,2}t_r\int_0^\infty g(t_ry) \{(s_r-u)^+ e^{iy\sgn_{z_r}  \tz_r}+1-\max(s_r,u)\} 
 \\
 \nonumber&&\qquad e^{iy\sgn_{z_r}(z_r-s_r \tz_r)}e^{ \rho_{z_r,\tz_r, s_r}(y)} \ud y
 \eeq
 
 The proof of the fact that under  analyticity hypotheses,  the limit points are uniquely prescribed by the above  equations goes along the same lines as the proofs in Section 5.2 of   \cite{ACFTCL}, sketched as follows. First, we have to consider separately the case where $\Phi$ satisfies \eqre{uniquenessassumption2} and the case where $\Phi$ satisfies \eqre{phias}. In the case 
where $\Phi$ satisfies \eqre{uniquenessassumption2}, the proof is   very similar to the proof of the corresponding case in  Section 5.2 of   \cite{ACFTCL} and to  the detailed proof of the uniqueness for Lemma \re{204131} of the present paper,  using \eqre{eq:AssumpTau} instead of \eqre{ConditionOng}. The case where $\Phi$ satisfies \eqre{phias} is a little more delicate. As in 
Lemma 5.1 of   \cite{ACFTCL}, one first needs to notice that considered as functions of $t,t',t''$, the limit points satisfy an H\"older bound, using essentially the facts that 
for  any $2\kappa\in (0,\alpha/2)$
\begin{equation}\label{2310131}\limsup_{n\ge 1} \mathbb E\Big[ \big(\sum_{i=1}^n |a_{1i}|^2\big)^{2\kappa}\Big]<\infty\,,\end{equation}  and that for any $\bet\in(\al/1,1]$, there exists a constant $c = c(\alpha, \beta)$ such that for any $x, y$ in $\bC^-$, 
\begin{equation}\label{fondin231013}
| x^{\frac \alpha 2}  - y^ {\frac \alpha 2}  | \leq c |x - y | ^ { \beta} \left(  | x | \wedge | y |  \right)^ { \frac  \alpha 2 - \beta }.\end{equation}
Then one has to interpret the equation satisfied by the limit points  as a fixed point equation for a strictly contracting function in a space of H\"older functions: the key argument, to prove that the function is contracting, is to use the estimates given in Lemma 5.7 of \cite{charles_alice}.
\epr

 This concludes the proof of Proposition 
\re{PropoCVST} and 
 it follows from this that the covariance of the process $H_{s,z}$ is given by
\be\la{591300h29}C(s,z;s',z'):=\mathbb E[H_{s,z}H_{s',z'}]=\int_0^1 \ud u (L_u(s,z;s',z')-L_u(z,s)L_u(z',s'))\ee
with the functions $L$ defined in Lemmas \ref{204131}  and \ref{defL2}.

%\tred{Doit-on laisser la remarque ci-dessous sur l'annulation de la covariance, sachant que l'on a maintenant une preuve compl√ãtement rigoureuse plus haut ?}
%
%Note that the above variance  $C(s,z,s,\bar z)$ does not vanish for $z\in\mathbb C\backslash \mathbb R$. Indeed,
%this would amount to see that for almost all $k$ the covariance of $f_k$ conditionally to $\mathcal F_{k-1}$ vanishes, that is that
%there exists a function $C_k\in \mathcal F_{k-1}$
%$$\f{\one_{k\le ns}-s_n+\sum_j\ba_k(j)^2\{G^{(k)}(P_s^{(k)}-s_n)G^{(k)}\}_{jj} }{z  -\sum_j\ba_k(j)^2 G^{(k)}_{jj}}=C_k^s(z)$$
%In other words, there exists another constant $\tilde C_k$ so that
%$$\sum_j\ba_k(j)^2[\{G^{(k)}(P_s^{(k)}-s_n)G^{(k)}\}_{jj}-C_k ^s(z) G^{(k)}_{jj}]= \tilde C_k^s$$
%which implies 
%$$\{G^{(k)}(P_s^{(k)}-s_n)G^{(k)}\}_{jj}-C_k^s (z)G^{(k)}_{jj}=0$$
%for all $j$.  We can next sum over $j$ to conclude by concentration inequalities that $C_k^s(z)$ is almost deterministic 
%and then taking expectations to deduce from the fact that the left hand side is centered that $C_k^s(z)$ must vanish. But this is not possible for $s\in (0,1)$
%  
%
%
%
%
%

 \section{Proof of Proposition \re{NonZeroCovariance}}\la{sec:NonZeroCovariance}
Let us now prove that the limit covariance of $(C_{s,\lam})$ is not identically zero (hence this is also the case for $(B_{s,t})$ by \eqre{rel31513}). 
Using Lemma \re{lemCtightness} and \eqre{305131}, one easily sees that $(C_{s,\lam}^n)$ is uniformly bounded in $L^4$. It follows that \be\la{tfpasl0}\Var(C_{s,\lam})=\E[(C_{s, \lam})^2]=\lim_{n\to\infty} \E[(C_{s, \lam}^n)^2].\ee Thus we shall prove that the limit of $\E[(C_{s, \lam}^n)^2]$ is not identically zero.
 
 \emph{\underline{Preliminary computation:}} For  $x_1, \ld, x_n\in \C$ \st $x_1+\cd+x_n=0$,   for any $0\le \ell\le n$, we have \be\la{47131}\sum_{i=1}^\ell x_i=\sum_{i=1}^n\al_i x_i\qquad\trm{ for } \al_i:=\beg{cases} 1-\f{\ell}{n}&\trm{ if $i\le \ell$,}\\
-\f{\ell}{n} &\trm{ if $i> \ell$,}\en{cases}\ee
and that $\al_1+\cd+\al_n=0$.
Note also that for  $(X_1, \ld, X_n)$ an exchangeable  random vector and $\al_1,\ld, \al_n\in \C$ \st $\al_1+\cd+\al_n=0$, 
  we have \be\la{47132}\E \sum_{i,i'}\al_i\al_{i'} X_iX_{i'}=\sum_i\al_i^2\E[X_1(X_1-X_2)].\ee
 It follows   from \eqre{47131} and \eqre{47132} that if the coordinates of an exchangeable random vector $(X_1, \ld, X_n)$    sum up  to zero, then for any  $0\le \ell\le n$,   \be\la{47133}\E \sum_{\substack{1\le i\le \ell\\ 1\le i'\le \ell}}X_iX_{i'}
=n\lf(\f{\ell}{n} -\f{\ell^2}{n^2}\ri)\E[X_1(X_1-X_2)].\ee

Let us now fix   $s\in (0,1)$ and $\lam\in \R$ and  apply our preliminary computation \eqre{47133} with $X_i=\sum_{j\ste \lam_j\le \lam}(|u_{ij}|^2-n^{-1})$ and $\ell=\lfl ns\rfl$.
For $s_n:={\lfl ns\rfl}/{n}$, we get  \be\la{tGC20:08}\var (C^n_{s,\lam})= \E[(C^n_{s,\lam})^2] = \lf(s_n-s_n^2\ri)\E[X_1(X_1-X_2)]\ee
Note also that as each $|u_{ij}|^2$ has expectation $n^{-1}$, $$
\E[X_1(X_1-X_2)]=\E[\sum_{\substack{j\ste \lam_j\le \lam\\ j'\ste \lam_{j'}\le \lam}}|u_{1j}|^2|u_{1j'}|^2-|u_{1j}|^2|u_{2j'}|^2]
$$
Moreover, by exchangeability of the rows of $U$ (which is true even conditionally to the $\lam_j$'s) and the fact that its columns have norm one,   for any $j,j'$, 
$$n(n-1)\E[\one_{\lam_j, \lam_{j'}\le \lam}|u_{1j}|^2|u_{2j'}|^2]+n\E[\one_{\lam_j, \lam_{j'}\le \lam}|u_{1j}|^2|u_{1j'}|^2]=1,$$ so that $$\E[X_1(X_1-X_2)]=O\lf(\ff{n}\ri)+\lf(1-\ff{n}\ri) \E[\sum_{j,j'\ste \lam_j,\lam_{j'}\le \lam}(|u_{1j}|^2|u_{1j'}|^2-n^{-2})].$$
By \eqre{tGC20:08}, we deduce that $$\E[(C^n_{s,\lam})^2] =O\lf(\ff{n}\ri)+\lf(1-\ff{n}\ri)\lf(s_n-s_n^2\ri)\E[\sum_{j,j'\ste \lam_j,\lam_{j'}\le \lam}(|u_{1j}|^2|u_{1j'}|^2-n^{-2})],$$
 so that \be\la{67136h54}\E[(C^n_{s,\lam})^2] =O\lf(\ff{n}\ri)+ \lf(s_n-s_n^2\ri)\E[\bigg(\sum_{j\ste \lam_j \le \lam}(|u_{1j}|^2 -n^{-1})\bigg)^2]\ee
%Note that $$\sum_{j\ste \lam_j\le \lam}|u_{1j}|^2=\mu_{n,\bfe_1}((-\infty, \lam])$$where  $\ds \mu_{n,\bfe_1}:=\sum_{j=1}^n|u_{1j}|^2\del_{\lam_j}$ is the empirical spectral law of $A$ according to the first vector $\bfe_1$ of the canonical basis, also defined by the fact that for any test function $f$, $$\int f(x)\ud \mu_{n,\bfe_1}(x)=(f(A))_{11}.$$
Moreover,  for $\mu_n$,  $\mu_{n,\bfe_1}$  the random \pro measures  introduced in  \eqre{2441411h31}, we have $$\sum_{j\ste \lam_j\le \lam}(|u_{1j}|^2-n^{-1})=(\mu_{n,\bfe_1}-\mu_n)((-\infty, \lam])$$
Hence by \eqre{67136h54}, \be\la{tfpasl}\E[(C^n_{s,\lam})^2] =O\lf(\ff{n}\ri)+ \lf(s_n-s_n^2\ri)\E[\{(\mu_{n,\bfe_1}-\mu_n)((-\infty, \lam])\}^2].
\ee

Let us now suppose that for a certain $s\in (0,1)$, we have $\Var(C_{s,\lam})=0$ for all $\lam\in \R$. To conclude the proof, we shall now exhibit a contradiction. By \eqre{tfpasl0} and \eqre{tfpasl}, we know that for all $\lam$, $\E[\{(\mu_{n,\bfe_1}-\mu_n)((-\infty, \lam])\}^2]\lto 0$ as $n\to\infty$, hence  $\E[|(\mu_{n,\bfe_1}-\mu_n)((-\infty, \lam])|]\lto 0$. As $\mu_{n,\bfe_1}$, $\mu_n$ are \pro measures, for any $\lam\in \R$, $$|(\mu_{n,\bfe_1}-\mu_n)((-\infty, \lam])|\le 2.$$ 
Thus  for any $z\in \C\bck\R$, by dominated convergence, as $n\to\infty$, $$\int_{\lam\in\R}\f{\E[|(\mu_{n,\bfe_1}-\mu_n)((-\infty, \lam])|]}{|z-\lam|^2}\ud\lam\lto 0.$$ We deduce the convergence in \pro $$ \int_{\lam\in \R}\f{(\mu_{n,\bfe_1}-\mu_n)((-\infty, \lam])}{(z-\lam)^2}\ud\lam \lto 0.$$  
But by \eqre{67136h15}, for any $z\in \C\bck\R$, with the notation $G(z):=(z-A)^{-1}$, \beq \int_{\lam\in \R}\f{(\mu_{n,\bfe_1}-\mu_n)((-\infty, \lam])}{(z-\lam)^2}\ud\lam&=&-\int_{\lam\in\R}\f{\ud(\mu_{n,\bfe_1}-\mu_n)(\lam)}{z-\lam}\\
&=&\ff{n}\Tr G(z)-G(z)_{11}.\eeq
We deduce the convergence in probability, for any fixed $z\in \C^+$, \be\la{77131} \ff{n}\Tr (G(z))-G(z)_{11} \lto 0.\ee  By \eqre{77132}, we deduce that $G(z)_{11}$ converges in \pro to the Stieltjes transform $G_{\mu_\Phi}(z)$ of the limit empirical spectral law $\mu_\Phi$ of $A$.
By the Schur complement formula (see \cite[Lem. 2.4.6]{agz}) and the asymptotic vanishing of non diagonal terms in the quadratic form (Lemma 7.7 of \cite{ACFTCL}), we deduce the convergence in probability    \be\la{87131}z-\sum_{j=2}^n|a_{1j}|^2G^{(1)}(z)_{jj}\lto 1/G_{\mu_\Phi}(z),\ee  where $A^{(1)}$ is  the matrix obtained after having removed the first row and the first column to $A$ and   $G^{(1)}(z):=(z-A^{(1)})^{-1}$.

It follows that the (implicitly depending on $n$) random variable $X=\sum_{j=2}^n|a_{1j}|^2G^{(1)}(z)_{jj}$ converges in \pro to a deterministic limit as $n\to\infty$. Let us show that this is not possible if $\Phi$ is not linear.

Let $\E_G$ denote the integration with respect to the randomness of the first row of $A$.
The random variable $X$ takes values in $\C^-$.  For any $t\ge 0$, by Lemma \re{2541400h15}, we have \be\la{2541414h}\E[e^{-itX}]=\prod_{j=2}^n\phi_n(tG^{(1)}(z)_{jj})=(1+o(1))\exp\{\ff{n-1}\sum_{j=2}^n\Phi(tG^{(1)}(z)_{jj})\}\ee
By Equation (20) of \cite{ACFTCL}, we know that \be\la{2541414h1}\ff{n-1}\sum_{j=2}^n\Phi(tG^{(1)}(z)_{jj})\lto \rho_z(t),\ee where $\rho_z$ is a continuous function on $\R^+$ satisfying (by Theorem 1.9 of \cite{ACFTCL}):
\be\la{2541414h2}\rho_z(\lam)=\lam\int_0^{+\infty}g(\lam y)e^{iyz+\rho_z(y)}\ud y.\ee
By \eqre{2541414h} and \eqre{2541414h1}, as $n\to\infty$, 
$$\E[e^{-itX}]\lto e^{ \rho_z(t)}.$$
but we already saw that $X$ converges in \pro to a constant, hence there is $c_z\in \C$ \st for all $t$, $\rho_z(t)=c_zt.$
From \eqre{2541414h2} and \eqre{hypcalcconv}, we deduce that for all $\lam\ge 0$, $$c_z\lam=\lam\int_0^{+\infty}g(\lam y)e^{iyz+c_zy}\ud y = \int_0^{+\infty}g(t)e^{i\f{z-ic_z}{\lam}t}\ud t=\Phi(\f{\lam}{z-ic_z}).$$
As we supposed that $\Phi$ is not linear, by analytic continuation, this is a contradiction. 

Note that the fact that $E[(C^n_{s,\lam})^2] $ does not go to zero could also be deduced from \cite{BCCAOP} in the \Lvy case according to \eqref{tfpasl}.

%\section{Proof of Proposition \re{791316h37}}\la{secProof791316h37}
%
%Recall that  $(C_{s,\lam})_{s\in [0,1], \lam\in \R}$ is the limit distribution of $(C^n_{s,\lam})_{s\in [0,1], \lam\in \R}$ as $n\to\infty$, hence for all 
%
%From the fact that $C$ is \cadlag in its second variable,   we have, for any $s,\lam$,  
%$$\tilde{C}_{s,\lambda}:=\f{1}{2}\bigg(\lim_{\tilde{\lam}\underset{<}{\to\lam}}C_{s,\tilde{\lam}}+C_{s, \lam}\bigg)=\lim_{\eta\downarrow 0}-\ff{\pi}\int_{-\infty}^\lam\Im(H_{s,E+i\eta})\ud E.$$
%We deduce that 
%$$\E[\tilde{C}_{s,\lambda}\tilde{C}_{s',\lambda'}]=\ff{\pi^2}\lim_{\eta\downarrow 0}
%\int_{-\infty}^\lambda\int_{-\infty}^{\lambda'}\mathbb E[ \Im \left( H_{s,E+i\eta}\right) \Im \left( H_{s',E'+i\eta}\right)]\ud E\ud E'\,.$$ 

\section{Appendix}
%\subsection{A lemma about exchangeable processes}
%First note that for  $x_1, \ld, x_n\in \C$ \st $x_1+\cd+x_n=0$,   for any $0\le \ell\le n$, we have \be\la{47131}\sum_{i=1}^\ell x_i=\sum_{i=1}^n\al_i x_i\qquad\trm{ for } \al_i:=\beg{cases} 1-\f{\ell}{n}&\trm{ if $i\le \ell$,}\\
%-\f{\ell}{n} &\trm{ if $i> \ell$,}\en{cases}\ee
%and that $\al_1+\cd+\al_n=0$.
%Note also that for  $(X_1, \ld, X_n)$ and $(Y_1, \ld, Y_n)$ be two   random vectors \st for any $i,j$, $\E[X_iY_j]$ depends only on weither $i$ equals $j$ or not %$\si\in S_n$, $$(X_{\si(1)},Y_{\si(1)},\ld\ld, X_{\si(n)},Y_{\si(n)})\eqlaw (X_{1},Y_{1},\ld\ld, X_{n},Y_{n})$$ 
%and for $\al_1, \ld, \al_n,\bet_1, \ld, \bet_n\in \C$ \st almost surely, $\bet_1+\cd+\bet_n=0$, 
%  we have \be\la{47132}\E \sum_{i,i'}\al_i\bet_{i'} X_iY_{i'}=\E[X_1(Y_1-Y_2)]\sum_i\al_i\bet_i.\ee
% It follows easily from \eqre{47131} and \eqre{47132} that if we also suppose that almost surely,  the $X_i$'s and the $Y_i$'s sum up (separately) to zero, then for any  $0\le \ell, \ell'\le n$, we have \be\la{47133}\E \sum_{\substack{1\le i\le \ell\\ 1\le i'\le \ell'}}X_iY_i
%=n\E[X_1(Y_1-Y_2)]\lf(\f{\ell}{n}\wedge \f{\ell'}{n}-\f{\ell}{n}\f{\ell'}{n}\ri).\ee

\subsection{A tightness lemma for   bivariate processes}
 Let us endow  the space $D([0,1]^2)$ with the Skorokhod topology (see \cite{FloEigenvectors} for the definitions). 

Let $M=[m_{ij}]_{1\le i,j\le n}$ be a random bistochastic matrix   depending implicitly on $n$. We define the random process of $D([0,1]^2)$ $$S^n(s,t):=\ff{\sqrt{n}}\sum_{\substack{1\le i\le ns\\ 1\le j\le nt}}\lf(m_{ij}-\ff{n}\ri).$$

\beg{lem}\la{lemCtightness}Let us suppose that $M$ is, in law, invariant under left multiplication by any permutation matrix. Then the process $S^n$ is $C$-tight in $D([0,1]^2)$, \ie  the sequence $(\op{distribution}(S^n))_{n\ge 1}$ is tight and has all its accumulation points   supported by the set of continuous functions on $[0,1]^2$. Moreover, the process $S^n$ is uniformly bounded in $L^4$.
\en{lem}

\bpr Let us prove that for all $0\le s<s'\le 1$, $0\le t<t'\le 1$,   \be\la{29513.8h} \E[( \Del_{s,s',t,t'} S^n)^4]\le \f{7}{n}+6(s'-s)^2(t'-t)^2(1-(s'-s))^2,\ee where $\Del_{s,s',t,t'} S^n$ denotes the increment of $S^n$ on $[s,s']\ti[t,t']$, \ie \be\la{KDdef29513}\Del_{s,s',t,t'} S^n:=\ff{\sqrt{n}}\sum_{\substack{ns< i\le ns'\\ nt\le j\le nt'}}(m_{ij}-\ff{n}).\ee
As $S^n$ vanishes on the boundary on $[0,1]^2$,    according to \cite[Th. 3]{bickel}, \eqre{29513.8h} will imply the lemma.

To prove \eqre{29513.8h}, we fix $0\le s<s'\le 1$, $0\le t<t'\le 1$. Let us now introduce some notation (where the dependence on $n$ will by dropped for readability). We define the sets $$I:=\{i\in \{1, \ld, n\}\ste ns< i\le ns'\} \trm{ and }J:=\{j\in \{1, \ld, n\}\ste nt< j\le nt'\},$$
the numbers $(\al_i)_{1\le i\le n}$ defined by $$\al_i:=\beg{cases}-\ff{\sqrt{n}}\f{|I|}{n}&\trm{ if $i\notin I$}\\ \\
\ff{\sqrt{n}}\lf(1-\f{|I|}{n}\ri)&\trm{ if $i\in I$}
\en{cases}$$and the exchangeable random vector (implicitly depending on $n$) $(X_1, \ld, X_n)$ defined by $$X_i=\sum_{j\in J}m_{ij}.$$
Note that\beq \Del_{s,s',t,t'} S^n&=&\ff{\sqrt{n}}\lf(\lf(\sum_{i\in I} X_i\ri)- \f{|I||J|}{n}\ri)
\eeq
and that as columns of $M$ sum up to one, $\ds |J|=\sum_{j\in J}\sum_{i=1}^nm_{ij}= \sum_{i=1}^nX_i$, 
hence $$
 \Del_{s,s',t,t'} S^n= \ff{\sqrt{n}}\lf( \sum_{i\in I} X_i- \f{|I|}{n}\sum_{i=1}^n X_i\ri)=\sum_{i=1}^n \al_i X_i.$$
 Thus by exchangeability of the $X_i$'s, we have \beq \E[( \Del_{s,s',t,t'} S^n)^4]
 &=&\E(X_1^4)\op{sum}_{4}(\al)+4\E(X_1^3X_2)\op{sum}_{3,1}(\al)+3\E(X_1^2X_2^2)\op{sum}_{2,2}(\al)\\&&+6\E(X_1^2X_2X_3)\op{sum}_{2,1,1}(\al)+\E(X_1X_2X_3X_4)\op{sum}_{1,1,1,1}(\al),\eeq with $$\op{sum}_{4}(\al):=\sum_{i=1}^n \al_i^4,\quad  \op{sum}_{3,1}(\al):=\sum_{i\ne j}\al_i^3\al_j,\quad  \op{sum}_{2,2}(\al):=\sum_{i\ne j}\al_i^2\al_j^2,$$ $$  \op{sum}_{2,1,1}(\al):=\sum_{\substack{i,j,k\\ \trm{pairwise $\ne$}}}\al_i^2\al_j\al_k,\quad   \op{sum}_{1,1,1,1}(\al):=\sum_{\substack{i,j,k,\ell\\ \trm{pairwise $\ne$}}}\al_i\al_j\al_k\al_\ell .$$
 As the $\al_i$'s sum up to zero, we have \beq \op{sum}_{3,1}(\al)&=&\sum_i(\al_i^3\sum_{j\ne i}\al_j)=-\op{sum}_{4}(\al),\\  \op{sum}_{2,1,1}(\al)&=&\sum_i(\al_i^2\sum_{j\ne i}(\al_j\sum_{k\notin \{i,j\}}\al_k)=\sum_i(\al_i^2\sum_{j\ne i}(\al_j(-\al_i-\al_j))\\&=&-\op{sum}_{3,1}(\al)-\op{sum}_{2,2}(\al)=\op{sum}_{4}(\al)-\op{sum}_{2,2}(\al)\\
 \op{sum}_{1,1,1,1}(\al)&=&-3\op{sum}_{2,1,1}(\al)=3\op{sum}_{2,2}(\al)-3\op{sum}_{4}(\al)
 \eeq
 Thus  \beq \E[( \Del_{s,s',t,t'} S^n)^4]
 &=&\op{sum}_{4}(\al)(\E(X_1^4)-4\E(X_1^3X_2)+6\E(X_1^2X_2X_3)-3\E(X_1X_2X_3X_4)) \\ &&+\op{sum}_{2,2}(\al)(3\E(X_1^2X_2^2)-6\E(X_1^2X_2X_3)+3\E(X_1X_2X_3X_4)),\eeq 
 Now, as for all $i$, $|\al_i|\le \ff{\sqrt{n}}$, we have $\op{sum}_{4}(\al)\le \ff{n}$, and as for all $i$, $0\le X_i\le 1$ (because the rows of $s$ sum up to one), we have \beq \E[( \Del_{s,s',t,t'} S^n)^4]
 &\le &\f{7}{n}  +3\op{sum}_{2,2}(\al)(\E(X_1^2X_2^2)+\E(X_1X_2X_3X_4)).\eeq 
 To conclude the proof of \eqre{29513.8h}, we shall prove that  
 \be\la{ihptalk1}  \op{sum}_{2,2}(\al)\le  (s'-s)^2  \ee
 and  \be\la{ihptalk2}  E(X_1^2X_2^2)+ \E(X_1X_2X_3X_4) \le 2 (t'-t)^2.\ee
 Let us first check \eqre{ihptalk1}. We have $$\op{sum}_{2,2}(\al)\le (\sum_i \al_i^2)^2= \lf\{(n-|I|)\f{|I|^2}{n^3}+|I|\ff{n}(1-|I|/n)^2\ri\}^2=\lf\{\f{|I|}{n}(1-\f{|I|}{n})\ri\}^2,$$ which gives \eqre{ihptalk1}.  Let us now check  \eqre{ihptalk2}. As $0\le X_i\le 1$, it suffices to prove that \be\la{ihptalk3}  \E(X_1X_2) \le  (t'-t)^2.\ee
 We have $$\E(X_1X_2)=\sum_{j,j'\in J} \E(m_{1j}m_{2j'}), $$ so it suffices to prove that   uniformly on $j,j'\in \{1, \ld, n\}$, $\ds \E(m_{1j}m_{2j'})\le \f{ 1}{n(n-1)}.$ We get this as follows: using the exchangeability of the rows of $M$ and the fact that  its rows sum up to one, we have, for any $j,j'\in \{1, \ld, n\}$, $$1=\E((\sum_i m_{ij})(\sum_{i'} m_{i'j'}))=n(n-1)\E(m_{1j}m_{2j'})+n\E(m_{1j}m_{1j'}).$$
 This concludes the proof.
\epr

\subsection{Injectivity of the Cauchy transform for certain classes of functions}
 
\beg{lem}\la{lemcauchytransfo}Let $f$ be a real valued bounded \cadlag function on $\R$ vanishing at infinity  with at most countably many discontinuity points. Then $f$ is entirely determined by the function $$\ds K_f(z):=\int \f{f(\lam)}{(z-\lam)^2} \ud \lam \qquad\trm{($z\in \C\bck \R$)}.$$   
More precisely, for any $\lam\in \R$, we have \be\la{891312h18}f(\lam) =\lim_{\substack{\tilde{\lam}\downarrow\lam\\ \trm{$f$ is cont. at $\tilde{\lam}$}}}\lim_{\eta\downarrow 0} \ff{\pi} \int_{-\infty}^{\tilde{\lam}}\Im K_f(E+i\eta)\ud E.\ee
\en{lem}
 
\bpr  Let us introduce the Cauchy transform of $f$, defined, on $\C\bck\R$, by $\ds H_f(z):=\int \f{f(\lam)}{ z-\lam } \ud \lam$.  It is well known that 
at any $\tilde{\lam} $ where $f$ is continuous, we have $$f(\tilde{\lam} )=\lim_{\eta\downarrow 0} -\ff{\pi}\Im H_f(\tilde{\lam} +i\eta).$$ Then, the result follows because  for all $\tilde{\lam}\in \R$, $\eta>0$, $$-H_f(\tilde{\lam}+i\eta)=  \int_{-\infty}^{\tilde{\lam}}K_f(E+i\eta)\ud E.$$ 
\epr 

%\beg{lem} Let $\nu$ be a complex measure on $\R$ with   Cauchy transform $$H_\nu(z):=\int \f{\ud \nu(t)}{z-t}\qquad \trm{ ($z\in \C\bck\R$)}.$$ Then for any $\lam\in \R$, 
%\be\la{79131}\f{\nu((-\infty, \lam))+\nu((-\infty, \lam])}{2}=\lim_{\eta\downarrow 0}-\ff{\pi}\int_{-\infty}^\lam \Im (H_\nu(u+i\eta))\ud u.\ee
%\en{lem}
%
%\bpr Let us fix $\lam\in \R$. As any complex measure is a linear combination of the Dirac mass $\del_\lam$ and of some of \pro measures  with no atom at $\lam$, it suffices to treat both of these case. In both of these cases, for any $\eta>0$, the function $u\in \R\longmapsto -\ff{\pi}  \Im (H_\nu(u+i\eta))$ is the density of the law of $X+\eta C$, where $X$ is a $\nu$-distributed random variable and $C$ is a random variable with Cauchy distribution, independent of $X$. Thus if $\nu=\del_\lam$, \eqre{79131} is obvious by symmetry of the Cauchy distribution, whereas if $\nu$ is a law with no aton at $\lam$, \eqre{79131} follows from the fact that $X+\eta C$ converges in distribution to $\nu$ as $\eta\downarrow 0$.
%\epr

\subsection{A lemma about large products and the exponential function}
The following lemma helps controlling the error terms in the proof of Proposition \re{PropoCVST} (in these cases, $M$ always has order one).
\beg{lem}\la{2541400h15} Let $u_i$, $i=1, \ld, n$, be some complex numbers and set 
 $$P:=\prod_{i=1}^n (1+\f{u_i}{n})\qquad\qquad S:=\ff{n}\sum_{i=1}^n u_i.$$ There is a universal constant $R>0$ \st for $M:=\max_i|u_i|$, $$  \f{M}{n}\le R\;\implies\; |P-e^S|\le \f{M^2}{n}e^{|S|+\f{M^2}{n}} .$$
\en{lem}

\bpr Let $L(z)$ be defined on $B(0,1)$ by $\log(1+z)=z+z^2L(z)$ and $R>0$ be \st on $B(0,R)$, $|L(z)|\le 1.$
If $\f{M}{n}\le R$, we have  $$P\ =\ \prod_i\exp\lf\{\f{u_i}{n}+\f{u_i^2}{n^2}L(\f{u_i}{n})\ri\}\
=\ e^S\exp\lf\{\sum_i\f{u_i^2}{n^2}L(\f{u_i}{n})\ri\},$$
which allows to conclude easily, as  for any $z$, $|e^{z}-1|\le |z|e^{|z|}.$ \epr

\subsection{CLT for martingales}

Let $(\mc{F}_k )_{k\ge 0}$ be a filtration \st $\mc{F}_0 =\{\emptyset,\Omega\}$ and let $(M_k )_{k\ge 0}$ be a square-integrable complex-valued martingale starting at zero with respect to this filtration. For $k\ge 1$, we define the random variables $$Y_k:=M_k-M_{k-1}\qquad v_k:=\E[|Y_k|^2\,|\,\mc{F}_{k-1} ]\qquad  \tau_k:=\E[Y_k^2\,|\,\mc{F}_{k-1} ]$$ and we also define $$v:=\sum_{k\ge 1}v_k\qquad  \tau:=\sum_{k\ge 1} \tau_k\qquad L(\eps):=\sum_{k\ge 1} \E[|Y_k|^2\one_{|Y_k|\ge \eps}].$$

Let now everything depend on a parameter $n$, so that $\mc{F}_k=\mc{F}_k^{\,n}, M_k=M_k^n, Y_k=Y_k^n,v=v^n,\tau=\tau^n,  L(\eps)=L^n(\eps), \ldots$

Then we have the following theorem. It is proved in the real case at \cite[Th. 35.12]{Billingsley}.  The complex  case can be deduced noticing that for $z\in \C$, $\Re(z)^2, \Im(z)^2$ and $\Re(z)\Im(z)$ are linear combinations of $z^2$, $\ovl{z}^2$, $|z|^2$.
\beg{Th}\la{thconvmart}Suppose that for  some  constants $v\ge  0, \tau\in \C$, we have the convergence in probability $$v^n\ninf v\qquad \tau^n\ninf \tau$$ and that for each $\eps>0$, $$L^n(\eps)\ninf 0.$$ Then we have the convergence in distribution $$M^n_n\ninf Z,$$ where $Z$ is a centered complex Gaussian variable \st $\E(|Z|^2)=v$ and $ \E(Z^2)=\tau$. 
\en{Th}

\subsection{Some linear algebra lemmas}
Let $\|\cdot\|_\infty$ denote the operator norm of matrices associated with the canonical Hermitian norm. 
\beg{lem}\la{formule_lin_alg} Let $A=[a_{ij}]$ be an $n\ti n$ Hermitian matrix, $z\in \C\bck\R$, $\ds G:=(z-A)^{-1}$,   $P$ be a diagonal matrix. For 
$1\le k\le n$ we denote by  $A^{(k)}$, $P^{(k)}$ be the matrices  with size $n-1$ obtained by removing the $k$-th row and the $k$-th column of $A$ and $P$ and set $G^{(k)}:=(z-A^{(k)})^{-1}$.  Then \be\la{104138h}
\Tr(P G)-\Tr( P^{(k)}G^{(k)})=\f{P_{kk}+\ba_k^*G^{(k)}P^{(k)}G^{(k)}\ba_k}{z-a_{kk}  -\ba_k^* G^{(k)}\ba_k},\ee with $\ba_k$ the $k$-th column of $A$ where the diagonal entry has been removed.
Moreover,  \be\la{161312h} |\Tr(P G)-\Tr( P^{(k)}G^{(k)})|\le \f{5\|P\|_\infty}{|\Im z|}
.\ee  
\en{lem}

\bpr $\bullet$ Let us first prove \eqre{104138h}. By linearity, one can suppose that $P$ has only one nonzero diagonal entry, say the $i$th one, equal to one.  Using the well known formula $$((z-A)^{-1})_{ii}-\one_{i\ne k}((z-A^{(k)})^{-1})_{ii}=\f{G_{ki}G_{ik}}{G_{kk}},$$we have  
\beq  \Tr(P G)-\Tr( P^{(k)}G^{(k)})
&=& ((z-A)^{-1})_{ii}-\one_{i\ne k}((z-A^{(k)})^{-1})_{ii}\\
&=&  \f{G_{ki}G_{ik}}{G_{kk}}\\
&=&\f{   ((z-A)^{-1}P (z-A)^{-1})_{kk} }{((z-A)^{-1})_{kk}}\\
&=&\f{{\pa_t}_{|_{t=0}} ((z-A-tP)^{-1})_{kk}  }{((z-A)^{-1})_{kk}}
\eeq
Let $\log$ denote the determination of the $\log$ on $\C\bck\R^-$ vanishing at one. Then 
\beq  \Tr(P G)-\Tr( P^{(k)}G^{(k)})&=&{\pa_t}_{|_{t=0}} \log \{((z-A-tP)^{-1})_{kk} \}\\
&=&{\pa_t}_{|_{t=0}} \log \ff{z-a_{kk}-t\one_{k=i} -\ba_k^* (z-X^{(k)} -tP^{(k)})^{-1}\ba_k}\\
&=&-{\pa_t}_{|_{t=0}} \log (z-a_{kk}-t\one_{k=i} -\ba_k^* (z-X^{(k)} -tP^{(k)})^{-1}\ba_k)\\
&=&\f{-{\pa_t}_{|_{t=0}} (z-a_{kk}-t\one_{k=i} -\ba_k^* (z-X^{(k)} -tP^{(k)})^{-1}\ba_k)}{(z-a_{kk}-t\one_{k=i} -\ba_k^* (z-X^{(k)} -tP^{(k)})^{-1}\ba_k)_{|_{t=0}}}\\
&=& \f{P_{kk}+\ba_k^*G^{(k)}PG^{(k)}\ba_k}{z-a_{kk}  -\ba_k^* G^{(k)}\ba_k}
\eeq

 $\bullet$ Let us now prove \eqre{161312h} (the proof does not use \eqre{104138h}). One can suppose that $k=1$. Let us introduce $$\tilde{A}:=\bbm a_{11} & 0&\cd&0\\
0&&&\\
\vdots & &A^{(1)} &\\
0&&&\ebm$$ and define $\tilde{G}$ and $\tilde{G}^{(1)}$ as $ {G}$ and $ {G}^{(1)}$ with $\tilde{A}$ instead of $A$. 
We have  $$|\Tr(P G)-\Tr( P^{(1)}G^{(1)})|\le |\Tr (P(G-\tilde{G}))|+
|\Tr(P \tilde{G})-\Tr( P^{(1)}\tilde{G}^{(1)})|+ |\Tr (P^{(1)}(G^{(1)}-\tilde{G}^{(1)}))|.$$ Let us treat the terms of the RHT separately. 

The third term is null because  $\tilde{A}^{(1)}=A^{(1)}$.  
We have $$|\Tr (P(G-\tilde{G}))|\le \|P(G-\tilde{G})\|_\infty\op{rank}(G-\tilde{G})$$ which is $\le \f{4\|P\|_\infty}{|\Im z|}$ by the resolvant formula. At last,   as $P$ is diagonal and the matrix $z-\tilde{A}$ can be inverted by blocs, we have  $$|\Tr(P \tilde{G})-\Tr( P^{(1)}\tilde{G}^{(1)})|=  |P_{11}\tilde{G}_{11}|\le \f{\|P\|_\infty}{|\Im z|}.$$
\epr


\begin{thebibliography}{10}
  \bibitem{agz} G.~Anderson, A.~Guionnet, O.~Zeitouni \emph{An Introduction to Random Matrices}. Cambridge studies in advanced mathematics, {\bf 118} (2009).
  \bibitem{abap} A. Auffinger, G.  Ben Arous,   S.
              P{\'e}ch{\'e}  \emph{Poisson convergence for the largest eigenvalues of heavy
              tailed random matrices},
{Ann. Inst. Henri Poincar\'e Probab. Stat.},
  {\bf 45}, {2009}, 589--610
\bibitem{bai-silver-book} Z.~D.~Bai, J.W.~Silverstein \emph{Spectral analysis of large dimensional random matrices}, Second Edition, Springer, New York, 2009.
\bibitem{BaoPanZhouUnivWigner} Z. Bao, G. Pan,   W. Zhou \emph{Universality for a Global Property of the Eigenvectors of Wigner Matrices}, arXiv 1211.2507
\bibitem{BeffaraDonatiRouault}  V. Beffara, C. Donati-Martin, A. Rouault \emph{Bridges and random truncations of random matrices}
Random Matrices: Theory and Appl. Vol. 03, No. 02.
  \bibitem{BDG} S. Belinschi, A. Dembo,  A. Guionnet   \emph{Spectral measure of heavy tailed band and covariance randommatrices}, {Comm. Math. Phys.}, {\bf 289},  {2009},  1023--1055.
  \bibitem{BAGheavytails}  G. Ben Arous and  A. Guionnet \emph{The spectrum of heavy tailed random matrices}. Comm. Math. Phys. {\bf 278} (2008), no. 3, 715--751.
   \bibitem{FloEigenvectors}  F. Benaych-Georges \emph{Eigenvectors of Wigner matrices: universality of global fluctuations}
Random Matrices Theory Appl. Vol. 1 (2012), no. 4, 23 pp. 
\bibitem{FT}  F. Benaych-Georges, T. Cabanal-Duvillard
\emph{Marchenko-Pastur Theorem and Bercovici-Pata bijections for heavy-tailed or localized vectors} ALEA, Lat. Am. J. Probab. Math. Stat.  Vol. 9  (2012),  no. 2, 685--715. 

   \bibitem{ACFTCL}  F. Benaych-Georges,  A. Guionnet, C. Male \emph{Central limit theorems for linear statistics of heavy tailed random matrices}.   {Comm. Math. Phys.}, Vol. 329 (2014), no. 2, 641--686. 
      \bibitem{bickel} P.J. Bickel, M.J. Wichura \emph{Convergence criteria for multiparameter stochastic processes and some applications}, Ann. Math. Statist., 42(5):1656--1670, 1971.
         \bibitem{Billingsley} P. Billingsley \emph{Probability and Measure}, Wiley, third edition.
              \bibitem{BCC2} C.~ Bordenave, P.~ Caputo,     D.~ Chafa{\"{\i}}
    \emph{Spectrum of non-{H}ermitian heavy tailed random matrices},
  {Comm. Math. Phys.},
 {\bf 307}, {2011}, 513--560.
\bibitem{BCCAOP} C.~ Bordenave, P.~ Caputo,     D.~ Chafa{\"{\i}}
    \emph{Spectrum of large random reversible Markov chains: heavy-tailed weights on the complete graph}. Ann. Probab. 39 (2011), no. 4, 1544--1590.
 \bibitem{charles_alice} C. Bordenave,  A. Guionnet \emph{Localization and delocalization of eigenvectors for heavy-tailed random matrices}, {Probab. Theory Related Fields} Vol. 157 (3-4), 885--953 (2013).
 \bibitem{CB} J.-P. Bouchaud, P. Cizeau, \emph{Theory of L\'evy matrices} Phys. Rev. E 50 (1994).
 
 \bibitem{chapuy} G. Chapuy \emph{Random permutations and their discrepancy process}, in 2007 Conference on Analysis of Algorithms, AofA 07, Discrete Math. Theor. Comput. Sci. Proc. AH (2007), pp. 415--426.
 
 
       \bibitem{cat-alain10} C. Donati-Martin, A. Rouault \emph{Truncations of Haar unitary matrices, traces and bivariate Brownian bridge}, Random Matrices: Theory and Application (RMTA) Vol. 01, No. 01.
       \bibitem{Uerdos} L. Erd{\H{o}}s, A. Knowles, H.T. Yau,
                 J. Yin, \emph{Spectral statistics of {E}rd{\H o}s-{R}\'enyi {G}raphs {II}:
              {E}igenvalue spacing and the extreme eigenvalues},
 {Comm. Math. Phys.}, {\bf 314},  {2012},
 {587--640}.
       \bibitem{U1} L.
       Erd{\H{o}}s, J. Ram{\'{\i}}rez, 
              B. Schlein,   H.T. Yau 
    {Universality of sine-kernel for {W}igner matrices with a small
              {G}aussian perturbation},
   \emph{Electron. J. Probab.},
 {\bf 15},
     {2010},
    no. 18, 526--603.
   \bibitem{Udeloc}L.
       Erd{\H{o}}s, 
              B. Schlein,     H.T Yau
    \emph{Semicircle law on short scales and delocalization of
              eigenvectors for {W}igner random matrices},
   {Ann. Probab.},
{\bf 37},
     {2009},
    815--852.
   \bibitem{Urigid}L.
       Erd{\H{o}}s, 
                H.T. Yau,   J. Yin  \emph{Rigidity of eigenvalues of generalized {W}igner matrices},
   {Adv. Math.}, {\bf 229},
      {2012},
     {3},
      {1435--1515}.		
\bibitem{U2}L.
       Erd{\H{o}}s, 
                H.T. Yau,   J. Yin \emph{Universality for generalized {W}igner matrices with
              {B}ernoulli distribution},
   {J. Comb.},
   {\bf 2},
      {2011}, {15--81}		
      \bibitem{ESY11} L.
       Erd{\H{o}}s, 
              B. Schlein,     H.T. Yau  \emph{Universality of random matrices and local relaxation flow},
 {Invent. Math.},
{\bf 185},  {2011},
   {75--119}		
\bibitem{U3} L.
       Erd{\H{o}}s, J. Ramirez, 
              B. Schlein, T. Tao, V. Vu,    H.T. Yau  
    \emph{Bulk universality for {W}igner {H}ermitian matrices with
              subexponential decay},
   {Math. Res. Lett.},
   {\bf 17},
     {\bf 2010},
   {667--674}		
\bibitem{U4}L.
       Erd{\H{o}}s, J. Ramirez, S. P{\'e}ch{\'e}
              B. Schlein,    H.T. Yau 
  \emph{Bulk universality for {W}igner matrices},
    {Comm. Pure Appl. Math.},
{\bf 63}, {2010}, {895--925}
\bibitem{KKP96}  A. M. Khorunzhy, B. A. Khoruzhenko, L.A. Pastur \emph{Asymptotic properties of large random
matrices with independent entries}, {J. Math. Phys.} 37 (1996) 5033--5060.  

\bibitem{KSV04} O. Khorunzhy, M. Shcherbina, V. Vengerovsky
\emph{Eigenvalue distribution of large weighted random graphs},
J. Math. Phys. 45 (2004), no. 4, 1648--1672. 
    

       \bibitem{MAL122}
{C. Male}
\emph{The limiting distributions of large heavy Wigner and arbitrary random matrices},
 {arXiv:1111.4662v3 preprint}.
 \bibitem{mehta}
 M. Mehta, \emph{Random matrices}, {Pure and Applied Mathematics (Amsterdam)}, {\bf 142},  {Third}, {Elsevier/Academic Press, Amsterdam}, 2004.
      
      \bibitem{tirozzi} M.~Shcherbina, B.~Tirozzi   \emph{Central limit theorem for fluctuations of linear eigenvalue
              statistics of large random graphs}, {J. Math. Phys.},
  {\bf 51},  {2010}, {023523, 20}
  
\bibitem{silverstein-AOP-eigenvectors} J. W. Silverstein \emph{Weak convergence of random functions defined by the eigenvectors of sample covariance matrices}.  Ann. Probab. 18 (1990), no. 3, 1174--1194. 
\bibitem{slanina} F. Slanina 
\emph{Localization of eigenvectors in random graphs},
{Eur. Phys. B}, {2012}, {85:361}.
\bibitem{sash} A. Soshnikov  \emph{ Poisson statistics for the largest eigenvalue of Wigner random matrices with heavy tails.} Electron. Comm. Probab. 9 (2004) 82--91.
 \bibitem{TV} T. Tao, V. Vu \emph{Random matrices: universality of local eigenvalue statistics},  Acta Math., 	
 206 (2011), 127--204.
\bibitem{Tao} T. Tao \emph{The asymptotic distribution of a single eigenvalue gap of a
              {W}igner matrix}, {Probab. Theory Related Fields}, {\bf 157},  {2013},
  81--106.
\bibitem{TW1} C.
Tracy,   H. Widom \emph{Level-spacing distributions and the Airy kernel.}
Comm. Math. Phys. {\bf 159}  (1994), no. 1, 151--174. 
 
  
\bibitem{ZAK}
  {I. Zakharevich}
\emph{A generalization of {W}igner's law},
{Comm. Math. Phys.},
 {268},
{2006},
{2},
 {403--414}.
\bibitem{wigner}  E.P. Wigner \emph{On the distribution of the roots of certain symmetric
              matrices},
  {Ann. Math.},
    {\bf 67},
     {1958},
    {325--327}.
\end{thebibliography}
\end{document}